\RequirePackage[reqno]{amsmath}
\documentclass[10pt,reqno]{amsart}

\usepackage{graphicx,graphics,caption, subcaption, enumerate,color}

\usepackage{amsmath,amssymb,amsfonts,bm,empheq}
\usepackage{cases}

\usepackage{comment}

\usepackage{soul}
 \usepackage{cancel}
 \usepackage[normalem]{ulem}
\usepackage{todonotes}

\usepackage{amssymb,amsmath,amsthm,amsfonts,amstext}

\usepackage{footnote}
\usepackage{epstopdf}

\usepackage{url}
\usepackage[colorlinks,bookmarksopen,bookmarksnumbered,citecolor=red,urlcolor=red]{hyperref}

\usepackage{lastpage}

\usepackage{enumitem}

\usepackage{multirow}
\usepackage{color}
\makesavenoteenv{tabular}

\newtheorem{theorem}{Theorem}

\newtheorem{lemma}[theorem]{Lemma}
\newtheorem{remark}{Remark}
\newtheorem{property}[theorem]{Property}
  \numberwithin{equation}{section}

\def\Real{\mathbb{R}}

\newcommand{\norm}[1]{\left\Vert#1\right\Vert}
\newcommand{\inpd}[2]{\left\langle #1, #2 \right\rangle}

\newcommand{\abs}[1]{\left\vert #1\right\vert}

\def\Hess{\mathbf{H}}
\newcommand{\wt}[1]{\widetilde{#1}}

\newcommand{\set}[1]{ \left\{ {#1}\right\}}

\newcommand{\argmin}{\operatornamewithlimits{argmin}}
\def\mcM{\mathcal{M}}
\def\mcU{\mathcal{U}}

\def\tspace{\mathcal{T}}
\def\vvTk{(v^{(k+1)}\otimes{v^{(k+1)}})}
\def\vvT{(v\otimes v)}
\def\kpa2_2{\frac{\kappa^2}{2}}

\def\sigm2_2{\frac{\sigma^2}{2}}
\def\nab_dot{\nabla\cdot}

\def\titphi{\tilde\phi}

\def\titF{\widetilde{F}}
\def\l{\text{l}}
\def\n{\text{n}}
\def\r{\textbf{r}}

\begin{document}

\begin{center}
{\bf Convex Splitting Method for  the  Calculation of \\
Transition States of Energy Functional}

\

\

Shuting Gu \footnote{email: shutinggu2-c@my.cityu.edu.hk.},
Xiang Zhou \footnote{email: xiang.zhou@cityu.edu.hk. The research of XZ  was supported by the grants from the Research Grants Council of the Hong Kong Special Administrative Region, China (Project No. CityU 11304314,    11304715
and 11337216).  }
\par
\
\par

Department of Mathematics
\par
City University of Hong Kong\par
Tat Chee Ave, Kowloon \par Hong Kong SAR

\end{center}

\section*{abstract}

     Among numerical methods for partial differential equations arising from steepest descent dynamics of energy functionals (e.g., Allen-Cahn and Cahn-Hilliard equations), the convex splitting method is well-known to maintain unconditional energy stability for a large time step size. In this work, we show how to use the convex splitting idea to find transition states, i.e., index-1 saddle points of the same energy functionals. Based on the iterative minimization formulation (IMF) for saddle points (SIAM J. Numer. Anal., vol. 53, p1786, 2015), we introduce the convex splitting method to minimize the auxiliary functional at each cycle of the IMF. We present   a general principle of constructing convex splitting forms for these auxiliary functionals and show how to avoid solving nonlinear equations. The new numerical scheme based on the convex splitting method allows for  large  time step sizes. The new methods are tested for the one dimensional Ginzburg-Landau energy functional in the search of the Allen-Cahn or Cahn-Hilliard types of transition states.  We provide the numerical results of transition states for the two dimensional Landau-Brazovskii energy functional for diblock copolymers.

\

{{\bf Keywords}: transition state, saddle point, convex splitting method, iterative minimization formulation
}
\

{{\bf  Mathematics Subject Classification (2010)}  Primary 65K05, Secondary  82B05	  }

\section{Introduction}
 For an energy functional,   both its local minimizers and its unstable saddle points have  important physical meanings for many problems in
physics, chemistry, biology and material sciences.
The local minimizers correspond to the stable configurations
in physical models,
and  they manifest  themselves  as   steady states of the gradient flow
driven by the energy.
 For the spatially extended systems,
these flows appear mathematically as the time-dependent  partial differential equations(PDEs). These PDEs   reflect the true physical  dynamics of
 relaxations and they   also
serve as a convenient  computational model to calculate
the stable steady states.
For instance, as the well-known phase separation and transition models,
 the Allen-Cahn (AC) \cite{AllenCahn} and Cahn-Hillian (CH) \cite{CahnHilliard} equations are the steepest descent dynamics of the Ginzburg-Landau energy functional under $L^2$ and $H^{-1}$ norms, respectively.
Besides the local minima, the other critical points on an energy surface
also play crucial roles for certain problems, such as the  energy-barrier activated processes which escape from local minima
by crossing   saddle points.
The infrequent hoppings between neighbouring local minima, although randomly,
occur in a quite certain style  of travelling  through  transition states.
These transition states,
as  the bottlenecks  on the pathways  of
 activated processes,  belong to a class of saddle points
 with index one, i.e.,  the unstable critical point
 whose  Hessian  has only one negative eigenvalue.

 In this paper, we are interested in how to find  these index-1 saddle points
for a given smooth energy functional.
The search for  transition states, or index-1 saddle points,
faces many challenges.
A large number of  numerical methods have been
proposed and developed to address these challenges. There are also  many  applications of  these   algorithms  in computational chemistry and material sciences.
 Refer to  \cite{Schlegel2003,ZHANG:hn} for   review  of this topic.
One class of numerical methods is to search the so-called minimum energy path (MEP).
The points along an MEP with locally maximal
energy value are then the index-1 saddle points.
These path-finding methods include   \cite{String2002, String2007}  and  \cite{NEB1998,CI-NEB2000}.
The other class of methods is to
evolve a single state
  on the potential energy surface.
The essential  question  is     how to define some dynamics on the energy surface
 to converge to index-1 saddle points without knowing
multiple local minima {\it  a prior}.
The intuitive idea of using the  softest direction (corresponding to the minimal eigenvalue of the Hessian) to invert the force component
along this min-mode direction proves very useful (\cite{HBK2005,OKHAJ2004}) and
it was proposed probably as early as in 1970s in \cite{Crippen1971,cerjan1981}.
 Many well-known algorithms and softwares such as the
  dimer method (\cite{Dimer1999,KS2008,DuSIAM2012}) or  the activation-relaxation techniques (\cite{ART1998, Mousseau2000, Cances2009})
are based on this min-mode-following idea.

 The underlying dynamics  of these min-mode-following algorithms
   has been rigorously formulated  and analyzed in  \cite{GAD2011}.
This dynamics, with the name ``gentlest ascent dynamics'' (GAD),
 simultaneously evolves both a position variable and a direction variable.
By analyzing the eigenvalue of the GAD,
  \cite{GAD2011}    proved the locally linear convergence   to  saddle point.
To accelerate the convergence rate,
a discrete iterative mapping, named the ``iterative minimization formulation'' (IMF), has been  proposed in   \cite{IMF2014}.
IMF has   three advantageous features:
(1) it has the quadratic convergence rate for non-degenerate saddle points;
(2) it turns the  problem of searching unstable saddle points  into a series of minimization subproblems;
(3) there is no restriction for numerical methods to solve
minimization subproblems in the IMF. The only important issue in practice
for better efficiency is
 to use the adaptive stopping rule
 in solving the subproblems.
 For a thorough discussion of the practical  algorithms based
 on the IMF, the readers can refer to  \cite{IMA2015}.

The advantages and flexibilities offered by the IMF
  immediately provide many   new opportunities to   explore the
existing methods  which were designed for searching local minima.
 The convex splitting method, originally proposed in \cite{eyre1998unconditionally},
successfully gives unconditionally energy stable schemes
to  ensure a large time step size.
The effectiveness of this method in resolving
   gradient dynamics as well as calculating the minimizers
 has been demonstrated by a vast number of  applications, for example,  the phase field model in \cite{eyre1998marching},  the phase field crystal model in  \cite{wise2009energy},  the thin film epitaxy model in \cite{JieShen2012},   the binary fluid surfactant model in \cite{gu2014energy}, as well as many others (\cite{JieShen2014, CWang_mPFC2011, SMWise2010, SMWise2009}).

 In this paper,  our motivation is to test the performance of
 the convex splitting method if this strategy is  used to
 locate the saddle point.  As we mentioned earlier,
 the IMF solves the saddle point   problem by solving a series of minimization subproblems,
 and these subproblems can be solved by  running the steepest decent dynamics.
 Since the outer iteration of the IMF (referred  as ``cycle'' in \cite{IMF2014,IMA2015}) is of quadratic convergence rate and it
usually only takes a few cycles in practice to reach the desired accuracy,
then one can expect that a better method for the subproblem may gain a better
speedup in   efficiency.
The  benefit of the convex splitting idea here
is that the subproblems can be solved by a   large time step size. Therefore
it takes less steps in each cycle to improve the overall efficiency of  locating the saddle point in the IMF.

 The idea of the convex splitting method is quite simple, but
there are two   important practical  issues   when   applied to specific problems.
The first   is   the construction of a convex splitting form
for a given energy functional.  In theory(\cite{eyre1998unconditionally}),  there  always
 exists  a convex splitting form for any continuous functional.
 The explicit decomposition has to be sought  for specific problems.
  The second
 is that  one should try best to construct
 a {\it linear}
  time-implicit term  in the convex spitting scheme
  since
this can avoid solving  a nonlinear system  at each  time step.

 The contributions in  our work of  applying the convex splitting method to
 saddle point search problems include the following:
 (1) for any given convex splitting form of the original energy functional,
 we  show how to   obtain the corresponding  convex splitting form
 of the auxiliary functional in the IMF. This means that
we design  an automatic  procedure
 from the traditional convex splitting method for local minimizers
 to the convex splitting method for saddle points;
  (2) we shall see later that the auxiliary functional in the IMF consists of  multiple terms involving the original energy functional.
 By adapting different convex splitting forms for
 different terms in the auxiliary functional, we can ensure the time-explicit discretization for nonlinear terms and obtain a linear system.
We demonstrate how to achieve this  by the examples of
Ginzburg-Landau energy functional  and the Landau-Brazovskii energy functional. The condition is  that one
need to know at least one  convex splitting form with the time-explicit nonlinear term
for the original energy functional.

The rest of the paper is organized as follows. In Section \ref{sec2},
we review the IMF for the saddle point search problem and the convex splitting method.
In Section \ref{sec3}, we construct  the convex splitting method for saddle point search  problems
and discuss numerical issues.
 Section \ref{sec4} presents the detailed numerical schemes and substantial numerical results for
  the   Ginzburg-Landau energy functional in the $L^2$ and $H^{-1}$ metric, respectively, subject to
   the Neumann or periodic boundary condition and the Landau-Brazovskii energy functional with periodic boundary condition. The conclusion is drawn in Section \ref{sec5}.

\section{Review}\label{sec2}
In this section, we review two foundations of our method, which were
born apparently from two different areas.

\subsection{The IMF for the saddle point search}\label{sec2-1}
 We first recall the iteration minimization formulation (IMF) in \cite{IMF2014}.  Let $\mcM $ be a Hilbert space equipped with the norm $\norm{\cdot}$
  and the inner product $\langle \cdot,\cdot \rangle$. Suppose that $V(x): \mcM \rightarrow \mathbb{R}$ is a sufficiently smooth potential function, then the IMF
  is the following iteration for the position variable $x$ and the direction variable $v$
\begin{numcases}{}
v^{(k+1)} = \argmin_{\|u\|=1} \inpd{u} {H(x^{(k)})u},\label{IMF_direction}\\
x^{(k+1)} = \argmin_y L(y;x^{(k)},v^{(k+1)}),\label{IMF_position}
\end{numcases}
where
$$    H(x^{(k)}) = ~ \nabla^2 V (x^{(k)}  ), $$
and
\begin{equation}\label{L_orig}
    \begin{split}
    L(y;x^{(k)}, v^{(k+1)}) = ~&(1-\alpha) V(y) + \alpha V\left(y- (v^{(k+1)} \otimes v^{(k+1)}) (y-x^{(k)})\right) \\
    ~ & - \beta V\left(x^{(k)} + (v^{(k+1)} \otimes v^{(k+1)})(y-x^{(k)})\right).
    \end{split}
 \end{equation}
 $\alpha$ and $\beta$ are two parameters and $\alpha+\beta>1$. Two special choices for $\alpha$ and $\beta$ are:
(i) $(\alpha, \beta) = (2,0),$ then $L(y;x,v) = -V(y) + 2 V(y - (v\otimes v) (y-x))$;
 (ii) $(\alpha, \beta) = (0,2),$ then $L(y;x,v) = V(y) - 2 V(x + (v\otimes v) (y-x))$.
So, the general form of $L$ in \eqref{L_orig} should be a linear combination of three terms involved in these two extreme cases
and furthermore the coefficient of the first  term $V(y)$ is   determined to be   $1-\alpha$ by examining  the Hessian $\nabla^2_y L(y=x; x, v)$.
 The main properties of the auxiliary objective function $L(y;x,v)$ when $\alpha + \beta >1$ are listed
 here for reference.
 \begin{theorem}[\cite{IMF2014}]
 \label{Th_IMF}~Suppose that $x^*$ is a (non-degenerate) index-1 saddle point of a $C^4-$function $V(x)$, i.e.,
its all eigenvalues are $\lambda_1<0<\lambda_2\leq \cdots $,
and the auxiliary function $L$ is defined by $(\ref{L_orig})$ with $\alpha + \beta >1$, then\\
 $(1)$ there exists a neighbourhood $\mcU$ of $x^*$ such that for any $x\in\mcU$, $L(y;x,v)$ is strictly convex in $y\in\mcU$ and thus has a unique minimum in $\mcU$;\\
 $(2)$ define the mapping $\Phi: x\in\mcU \rightarrow \Phi(x)\in\mcU$ to be the unique minimizer of $L$ in $\mcU$ for any $x\in\mcU$.  Then the mapping $\Phi$ has only one fixed point $x^*$;\\
 $(3)$ the mapping $x\to \Phi(x)$ has a quadratic convergence rate.
 \end{theorem}

The IMF includes two-level iterations. The top level is   $x\rightarrow \Phi(x)$,
  referred   as  ``cycle''.
The $k$-th cycle means the step of $x^{(k)}\,\rightarrow \,x^{(k+1)}=\Phi(x^{(k)})$,
which in practice consists of a second-level iterative procedure to solve    \eqref{IMF_direction}
for the min-mode  (the so-called ``rotation step'') and
\eqref{IMF_position} to update the position(``translation step'').
The rotation step is a classical numerical eigenvector problem, for which
many methods have been constructed such as the power method (\cite{GAD2011}), the conjugate gradient method(\cite{PhysRevLett.86.664, HBK2005}),
the Lanczos algorithm(\cite{ART1998, Mousseau2000, Cances2009})
and the LOR in \cite{LOR2013}.

We are interested in   spatially extended systems, i.e., $\mcM$ is a function space and $V$ is actually a  functional on $\mcM$.  As a convention, we use $F(\phi)$ rather than $V(x)$ below to represent
the functional of a spatial function $\phi$.

\subsection{Convex splitting method}\label{sec2-2}
Let  $\phi(x,t): [0,1] \times \Real^+\to \Real$ be the solution of  the following PDE driven by the gradient flow
\begin{equation}\label{pde_ex}
\frac{\partial \phi}{\partial t} = -\frac{\delta F}{\delta \phi}(\phi),
\end{equation}
subject to certain boundary condition at $x=0$ and $1$.
 $F$ is a sufficiently smooth free energy functional bounded from below.
 $\frac{\delta F}{\delta \phi}$ is the first order variational derivative of $F(\phi)$ with respect to $\phi$.
 A convex splitting form of   $F(\phi)$  means  that
two convex functionals exist, denoted by $F_c$ and $F_e$,  such that
 \begin{equation}\label{convex_ex}
 F(\phi) = F_c(\phi) - F_e(\phi),
 \end{equation}
  where ``$c$'' (``$e$'') refers to the contractive (expansive) part of the energy (\cite{eyre1998unconditionally}). Then the convex splitting scheme for (\ref{pde_ex}) is
\begin{equation}\label{scheme_ex}
 \frac{\phi^{n+1} - \phi^n}{\Delta t} = -\left(
 \frac{\delta F_c}{\delta \phi}\left(\phi^{n+1}\right) - \frac{\delta F_e}{\delta \phi}
 \left(\phi^n\right)
  \right),
\end{equation}
where $\phi^n \approx \phi(t_n)$ is the numerical solution at the $n$-th time level  $t_n = n\Delta t$
and $\Delta t$ is the time step size.
 The time-discrete scheme (\ref{scheme_ex}) has the property of the so-called unconditional energy stability,
as stated in  the following theorem.

\begin{theorem}[\cite{wise2009energy}]
\label{thm:wise}~Suppose the free energy functional $F(\phi)$ can be split into two parts $F(\phi) = F_c(\phi)-F_e(\phi)$  as in $(\ref{convex_ex})$. Then the time-discrete scheme $(\ref{scheme_ex})$ is unconditionally energy stable, meaning that for any time step size $\Delta t > 0$,
$$ F(\phi^{n+1}) \leq F(\phi^n), \quad n = 0, 1, 2,\cdots. $$
\end{theorem}
The convexity of $F_c$ and $F_e$ do not have to be valid in the whole configuration space. In practical applications,
it only requires the convexity in bounded subsets of the configuration space, where the solution is known to be in.

  \section{Main method}\label{sec3}
Consider the  1-D spatial domain $(a,b)$.
The high dimensional case  follows exactly the same idea.
 Let  $\mathcal{M}$ be a function space on the interval $[a,b]$.
   For example, $\mcM$ is the Hilbert space $H^{1}([a,b])$ or
  some other subspaces of $L^2([a,b])$.
   { Our goal is to find the transition states of the energy functional $F$ on $\mcM$. }

  Assume that the second order  variational derivative (Hessian) of $F$ in $\mcM$, denoted by $H(\phi) := \frac{\delta^2 F(\phi)}{\delta\phi^2}$,  exists.
  We   rewrite the IMF in Section \ref{sec2-1}:
 \begin{numcases}{}
v^{(k+1)} = \operatornamewithlimits{argmin}_{\|v\|=1} \left\langle v, H(\phi^{(k)})v \right\rangle,\label{IMF_p1}\\
\phi^{(k+1)} = \operatornamewithlimits{argmin}_\phi L(\phi;\phi^{(k)},v^{(k+1)}),\label{IMF_p2}
\end{numcases}
 where     the auxiliary functional $L$ is
 \begin{equation}\label{L}
     L(\phi;\phi^{(k)}, v^{(k+1)}) = F(\phi)  -\alpha F(\phi)   + \alpha F (\phi- \titphi  ) \\
      - \beta F (\phi^{(k)}+ \titphi  ),
 \end{equation}
where we define
 \begin{equation}\label{eqn:tphi}
 \titphi := (v^{(k+1)}\otimes{v^{(k+1)}})(\phi-\phi^{(k)})
 \end{equation}
 to ease the notation.
  The min-mode $v^{(k+1)}$ belongs to $\tspace_{\mathcal{M}}$, the tangent space of $\mathcal{M}$. $u \otimes v$ denotes the tensor defined by $(u\otimes v)\phi =\left \langle v, \phi \right\rangle u$ for $u,v\in \tspace_{\mathcal{M}}$.

\subsection{Convex splitting method for  minimizing the auxiliary functional $L$ }\label{sec3-1}

We shall approximate the solution of the variational subproblem (\ref{IMF_p2})
$$ \min\limits_{\phi} L\left(\phi;\,\phi^{(k)},v^{(k+1)}\right) $$ at  the $k$-th IMF cycle
by the steady    solution of the gradient flow associated with  $L$,
 \begin{equation}\label{dynamics}
\frac{\partial\phi}{\partial t} = -\frac{\delta L}{\delta \phi}.
\end{equation}
Note $\phi^{(k)}$ and $v^{(k+1)}$ are  fixed here.  The solution  at  infinite time  is well-defined if the lowest eigenvalue of $L$ at $\phi = \phi^{(k)}$ is negative (\cite{IMF2014,IMA2015}).

Next, we show how to construct  the convex splitting scheme for equation (\ref{dynamics}).
Our starting point is that
a convex splitting form for $F(\phi)$ has been given, say,
   \begin{equation}\label{F_split}
F(\phi) = F_c(\phi)-F_e(\phi).
\end{equation}
  By substituting (\ref{F_split}) into (\ref{L}), we find the following convex splitting form for $L$  \begin{equation}\label{L_split}
 L(\phi;\, \phi^{(k)},v^{(k+1)})=L_c(\phi;\, \phi^{(k)},v^{(k+1)})-L_e(\phi;\, \phi^{(k)},v^{(k+1)}),
 \end{equation}
where  \begin{equation}\label{alpha_beta_geq01}
 \begin{split}
 & L_c(\phi;\phi^{(k)},v^{(k+1)}) =\\
  &\begin{cases}
    F_c(\phi)+ \alpha F_e(\phi) + \alpha F_c( \phi-\titphi ) + \beta F_e(\phi^{(k)} + \titphi ),  {\rm if} ~ \alpha \geq 0, \beta \geq 0; \\
    F_c(\phi)+ \alpha F_e(\phi) + \alpha F_c( \phi-\titphi ) - \beta F_c(\phi^{(k)} + \titphi),  {\rm if} ~ \alpha \geq 0,\beta \leq 0;\\
    F_c(\phi)- \alpha F_c(\phi) - \alpha F_e( \phi-\titphi ) + \beta F_e(\phi^{(k)} + \titphi),  {\rm if} ~ \alpha \leq 0,\beta \geq 0,
  \end{cases}
  \end{split}
  \end{equation}
and $ L_e(\phi;\phi^{(k)},v^{(k+1)})$ is defined likewise by exchanging the subindices
``$e$'' and ``$c$''  in \eqref{alpha_beta_geq01}.
Note that  $\titphi$ has   been defined in  \eqref{eqn:tphi}.


 \begin{property}\label{exist_convex_splitting}
 $L_c(\phi; \, \phi^{(k)},v^{(k+1)})$ and $L_e(\phi; \, \phi^{(k)},v^{(k+1)})$ defined above  are all convex with respect to $\phi$ for any $\phi^{(k)}$ and any $v^{(k+1)}$.
\end{property}
\begin{proof}
We only prove the case when $\alpha\geq 0, \beta\geq 0$ since the other two cases can be proved similarly.
The proofs of the convexity for  $L_c$ and $L_e$  are   the same, thus it  suffices to only show that two  special terms $F_c(\phi-\titphi)$ and $F_e(\phi^{(k)}+\titphi)$
in \eqref{alpha_beta_geq01} are both convex
in terms of $\phi$. Denote the Hessian of $F_c(\phi)$
and $F_e(\phi)$ by $H_c$ and $H_e$, respectively, then both $H_c$ and $H_e$ are semi-positive definite
by the property of convexity.
By the definition of  $\titphi$  in \eqref{eqn:tphi}, the second order derivatives of
$F_c(\phi-\titphi)$ and $F_e(\phi^{(k)}+\titphi)$ are
$
 [I-\vvT]H_c(\phi-\titphi)[I-\vvT] $
and
$ \vvT H_e(\phi^{(k)}+\titphi)\vvT,$
respectively, where $v=v^{(k+1)}$. Since $v$ is nonzero,
these two (projected) Hessians are also  semi-positive definite.
This completes our proof.
\end{proof}

We now present our time-discrete numerical scheme for \eqref{dynamics} based on the convex splitting idea
 in Section \ref{sec2-2}.
   Here we only consider the case $\alpha\geq0,  \beta \geq 0$.
To simplify the notation, we drop out  the parameters $\phi^{(k)}$ and $v^{(k+1)}$ in
the expressions of $L(\phi,\phi^{(k)},v^{(k+1)})$, $ L_c(\phi,\phi^{(k)},v^{(k+1)})$ and $L_e(\phi,\phi^{(k)},v^{(k+1)})$. Let $\phi^n$ be the numerical solution  at the time level $t_n$. Our  scheme is
\begin{equation}\label{time_discrete_scheme_1}
\frac{\phi^{n+1} - \phi^n}{\Delta t} = -\left[ \frac{\delta L_c}{\delta \phi}(\phi^{n+1}) - \frac{\delta L_e}{\delta \phi}(\phi^n) \right],
\end{equation}
with initial  $\phi^0 = \phi^{(k)}$.
The first order variational derivative  of $L_c$ is
 \begin{equation*}
    \begin{split}
        \frac{\delta L_c}{\delta \phi}(\phi) \!= & \frac{\delta F_c}{\delta \phi}(\phi) \!+ \!\alpha \frac{\delta F_e}{\delta\phi}(\phi)\!
        \\
         & + \alpha \big(I-\vvTk \big)\frac{\delta F_c}{\delta\phi}\left(\phi-\vvTk(\phi-\phi^{(k)})\right)
        \\& + \!\beta \vvTk \frac{\delta F_e}{\delta \phi} \left (\phi^{(k)}\!\! + \!\! \vvTk(\phi\!-\!\phi^{(k)}) \right ),
          \end{split}
\end{equation*}
and $\frac{\delta L_e}{\delta \phi}(\phi)$ has the similar form
by switching subindices ``$c$'' and ``$e$''.

For the same reason in  Theorem \ref{thm:wise},
  the following unconditional energy stability holds for the scheme (\ref{time_discrete_scheme_1}).
The proofs are straightforward and   skipped.

\begin{lemma}\label{inequality_lemma}
Suppose that $\phi, \ \psi: [a,b] \times \Real^+ \rightarrow \mathbb{R}$ are two  periodic functions. If $L(\phi)$ in $(\ref{L})$ has the convex splitting form $L(\phi) = L_c(\phi)-L_e(\phi)$ given in
\eqref{L_split}. Then
\[ L(\phi) - L(\psi)
\leq \inpd{ \delta_\phi L_c(\phi)- \delta_\phi L_e(\psi)}{\phi-\psi }_{L^2}, \]
where $\delta_\phi L_c$ and $\delta_\phi L_e $ represent the first order variational derivatives of $L_c$ and $L_e$ with respect to $\phi$, respectively.
\end{lemma}
\begin{theorem}
\label{egy_stable_thm}
 If the energy functional $F(\phi)$ has the convex splitting form $F=F_c-F_e$, then the time-discrete scheme $(\ref{time_discrete_scheme_1})$ is unconditionally energy stable, meaning that for any time step size $\Delta t >0$, we have
$$ L(\phi^{n+1}) \leq L(\phi^{n}), \quad n= 0,1, 2,\cdots, $$
in each $k$-th cycle.
\end{theorem}

\subsection{Avoid the time-implicit nonlinear term and construct the linear system }\label{sec3-2}
In many cases,  several convex splitting forms for   $F$
may be found.
 Assume $F(\phi)$ has two forms of convex splitting
 with the following properties:
  \begin{equation}
 F(\phi) = F_c^\text{l}(\phi) - F_e^\text{n}(\phi), \label{F_split_ln}
\end{equation}
and
 \begin{equation}
  F(\phi) = \widetilde F_c^\text{n}(\phi) - \widetilde F_e^\text{l}(\phi), \label{F_split_nl}
\end{equation}
where the superscripts ``l''
 and ``n''  mean that the first order variational derivative is linear or nonlinear in $\phi$, respectively.
Accordingly, two convex splitting schemes exist to solve the gradient flow   $\frac{\partial \phi}{\partial t} = - \frac{\delta F}{\delta \phi}$:
\begin{equation}\label{eqn:1f}
\frac{\phi^{n+1} - \phi^n}{\Delta t} = -
  \frac{\delta F^\text{l}_c}{\delta \phi} (\phi^{n+1}) +  \frac{\delta F^\text{n}_e}{\delta \phi}  (\phi^n),
\end{equation}
and
\begin{equation}\label{eqn:2f}
\frac{\phi^{n+1} - \phi^n}{\Delta t} =  - \frac{\delta \widetilde F^\text{n}_c}{\delta \phi} (\phi^{n+1}) +  \frac{\delta \widetilde F^\text{l}_e}{\delta \phi} (\phi^n).
\end{equation}
Both schemes satisfy the unconditional energy stability.
The difference
is that   the nonlinear terms are handled differently.
The scheme \eqref{eqn:1f}  is  time-explicit in  nonlinear term
and hence solves a linear system   to generate  $\phi^{n+1}$ at the next time level.
But the scheme \eqref{eqn:2f} requires to solve a nonlinear equation for $\phi^{n+1}$.
If one has  a very efficient nonlinear solver,  such as
the multi-grid method in \cite{Kornhuber2006} for vector-valued Allen-Cahn equation,
then  the scheme \eqref{eqn:2f} is also a good choice.
In other cases,  the linear scheme \eqref{eqn:1f}   is usually preferred.

In the context of the IMF, our proposed convex splitting scheme
for the saddle point problem faces the same difficulty of
possible emergence of  nonlinear time-implicit term.
If only  one convex splitting form
is available, then the convex splitting form of the auxiliary functional $L$
specified  by \eqref{alpha_beta_geq01}
 inevitably runs into this trouble.
Either $L_c$    or $L_e$    includes  both
$F_c$ and $F_e$, which  leads to the appearance of at least one nonlinear term
at the implicit time level.
However, if one has both  \eqref{F_split_ln} and  \eqref{F_split_nl},
then  this difficulty can be circumvented by combining them together.
Take $\alpha, \beta \geq 0$ as an example  again.
   Substituting \eqref{F_split_ln} into the first and the third terms on the right hand side of \eqref{L} and substituting \eqref{F_split_nl} into the second and the fourth terms on the right hand side of (\ref{L}), then
   we have the following decomposition of $L$:
\begin{align*}
L(\phi) = & ~~F_c^\text{l}(\phi) - F_e^\text{n}(\phi) - \alpha \big[ \widetilde F_c^\text{n}(\phi) - \widetilde F_e^\text{l}(\phi) \big] \\
 & + \alpha \big[ F_c^\text{l}(\phi-\titphi) - F_e^\text{n}(\phi-\titphi) \big] - \beta \big[ \widetilde F_c^\text{n}(\phi^{(k)} + \titphi) - \widetilde F_e^\text{l}(\phi^{(k)} + \titphi) \big]\\
 = &~\Big[ F_c^\text{l}(\phi) + \alpha\widetilde F_e^\text{l}(\phi) + \alpha F_c^\text{l}(\phi-\titphi) + \beta\widetilde F_e^\text{l}(\phi^{(k)} + \titphi) \Big] \\
 & - \Big[ F_e^\text{n}(\phi) + \alpha \widetilde F_c^\text{n}(\phi) + \alpha F_e^\text{n}(\phi-\titphi) + \beta \widetilde F_c^\text{n}(\phi^{(k)} + \titphi) \Big]\\
 =: &~ L_c(\phi) - L_e(\phi),
\end{align*}
where $\titphi = \vvTk (\phi - \phi^{(k)}).$ It is easy to see that both $L_c$ and $L_e$ are convex with respect to $\phi$. The proof is exactly the same as that for Property \ref{exist_convex_splitting}. Theorem \ref{egy_stable_thm} in Section \ref{sec3-1} also holds true in this case.

In constructing two convex splitting forms satisfying the above conditions,
we do not require  the convexity holds for the whole configuration space.
In fact, for many examples,  the  global convexity for  all
four functions
in \eqref{F_split_ln} and \eqref{F_split_nl}
is not possible.
 We only   need  the   convexity properties valid  locally  for a certain range of the   values of $\phi$
in which  the solution of interest does not violate.

In summary, the semi-discrete scheme for (\ref{dynamics}) by using the two forms of convex splitting of $F(\phi)$, \eqref{F_split_ln}  and \eqref{F_split_nl}, has the following expression:
\begin{equation}\label{time_discrete_scheme_2}
\frac{\phi^{n+1} - \phi^n}{\Delta t} = -\left[ \frac{\delta L_c}{\delta \phi}\right]^{n+1} + \left[ \frac{\delta L_e}{\delta \phi} \right]^n,
\end{equation}
where
\begin{equation*}
    \begin{split}
        \left[\frac{\delta L_c}{\delta \phi}\right]^{n+1} = & ~ \frac{\delta F_c^\l}{\delta \phi}(\phi^{n+1}) + \alpha \frac{\delta \widetilde F_e^\l}{\delta\phi}(\phi^{n+1})\\
         & + \beta \vvTk \frac{\delta \widetilde F_e^\l}{\delta \phi}
         \left(\phi^{(k)} + \vvTk(\phi^{n+1}\!-\!\phi^{(k)})\right)\\
         & + \alpha (I-\vvTk)\frac{\delta F_c^\l}{\delta\phi}
         \left(\phi^{n+1}-\vvTk(\phi^{n+1}-\phi_k)\right),\\
         \left[\frac{\delta L_e}{\delta \phi}\right]^n = & ~ \frac{\delta F_e^\n}{\delta \phi}(\phi^n) + \alpha \frac{\delta \widetilde F_c^\n}{\delta\phi}(\phi^n) \\
         & + \beta \vvTk \frac{\delta \widetilde F_c^\n}{\delta \phi}
         \left(\phi^{(k)} + \vvTk(\phi^n - \phi^{(k)})\right)\\
  & + \alpha (I-\vvTk)\frac{\delta F_e^\n}{\delta\phi}\left(\phi^n - \vvTk (\phi^n-\phi^{(k)})\right), \\
    \end{split}
\end{equation*}
 and the initial value $ \phi^0 = \phi^{(k)}.$ $\alpha\geq 0, \beta \geq 0$ and $\alpha+\beta>1$.
It is easy to see that $\Big[\frac{\delta L_c}{\delta \phi}\Big]^{n+1}$ is linear and $\Big[\frac{\delta L_e}{\delta \phi}\Big]^n$ is nonlinear. So the scheme (\ref{time_discrete_scheme_2}) indeed corresponds to a linear system.

\section{Applications}\label{sec4}
\subsection{ 1D example: Ginzburg-Landau free energy}\label{sec4_1}
We apply our method to  the   Ginzburg-Landau free energy on  $[0,1]$,
 \begin{equation}
 \label{eqn:F}
F(\phi) = \int_0^1 \Big[ \frac{\kappa^2}{2}(\frac{\partial \phi}{\partial x})^2 + f(\phi) \Big]\,dx,
\end{equation}
 where $\phi(x)$ is an order parameter representing for example the concentration
of one of the component in a binary alloy.
 The mobility parameter $\kappa>0$. $f(\phi) = (\phi^2-1)^2/4$.
If we consider   the gradient flow of $F$ in the  $L^2([0,1])$ space  with the standard $L^2$ inner product
$\inpd{\cdot\,}{\,\cdot}_{L^2}$,
then we obtain the (non-conserved) Allen-Cahn  (AC) equation
\begin{equation}\label{AC_eq}
\frac{\partial\phi}{\partial t} = -\frac{\delta F}{\delta\phi}(\phi)=\kappa^2 \Delta \phi - f'(\phi)
=\kappa^2 \Delta \phi - (\phi^3-\phi),
\end{equation}
where $\Delta = \partial_{xx}$.
If the gradient flow is defined  in the $H^{-1}$ metric
$\inpd{\cdot\,}{\,\cdot}_{H^{-1}}$,
then we have the (conserved) Cahn-Hilliard (CH) equation
\begin{equation}\label{CH_eq}
\frac{\partial\phi}{\partial t} = \Delta\frac{\delta F}{\delta\phi}=-\kappa^2 \Delta^2 \phi + \Delta (\phi^3-\phi).
\end{equation}

We are interested in the unstable index-1 saddle point of the Ginzburg-Landau free energy
\eqref{eqn:F}.  These saddle points correspond to the ``spike-like'' stationary solutions,
or ``canonical nuclei'' discussed in \cite{BatesFife1993}.
 Similarly  to the AC and CH equations, which arise  in $L^2$ metric and $H^{-1}$ metric, respectively,
we search for the saddle points   of $F$  both in $L^2$ and in $H^{-1}$ metrics.
The calculations of transition states and transition rates for the CH equation have already been done in \cite{CH1D2012,MMSLZZ2013}
by using the string method
(\cite{String2002}) and
the GAD(\cite{GAD2011}).

We consider  both  the Neumann and periodic   boundary condition.
The Neumann boundary condition  is $\partial_x\phi(0)=\partial_x\phi(1)=0$
for the AC equation and $\partial_x\phi(0)=\partial_x\phi(1)=\partial^3_x \phi(0)=\partial^3_x\phi(1)=0$
for the CH equation.
The periodic boundary condition simply means that $\phi(x)=\phi(x+1),~\forall x\in[0,1]$, which induces a degeneracy at
any stationary  solution corresponding to the    translation invariance    in the spatial variable $\phi(x)\to \phi(x+c)$.
 This means that the second smallest eigenvalue
  of the index-1 saddle points is zero.
The degeneracy from the periodic boundary condition does not affect the quadratic convergence rate of the IMF
(\cite{GUthesis}).
It is also noted that the mass $\int_0^1 \phi(x)dx$ is conservative in $H^{-1}$ metric and
  any stationary solution is still stationary if
 an arbitrary constant is added.
 This degeneracy
  can be    eliminated     by  restricting the   solutions
in the space where  the    mass  is chosen beforehand; for the same reason,
 any eigenvectors or perturbations should   be restricted  to  having zero  mass.
We refer the reader to  \cite{BatesFife1990, BatesFife1993}
and references therein on the existence of index-1  saddle point
for sufficient small $\kappa$.
We restrict  our calculation to the case of not too large domain,
i.e., the parameter $\kappa$ in \eqref{eqn:F} is not too small, but small enough
to posses saddle points.

Now we discuss some details of  our   method.
In the IMF, the auxiliary functional
$L$   given by \eqref{L} (setting $\alpha=0$, $\beta=2$)
is \begin{equation}\label{L-ACH}
L(\phi; \phi^{(k)}, v^{(k+1)})
=F(\phi)- 2 F(\hat\phi)
\end{equation}
with
\begin{equation}\label{370}
\hat{\phi}:=\phi^{(k)}+\titphi =\phi^{(k)} +  \vvTk (\phi - \phi^{(k)}),
\end{equation}
where   $(v\otimes v)u = \langle v,u \rangle v, \forall u, v,$ is associated with either $L^2$ metric or $H^{-1}$ metric. So the formal notation $\inpd{\cdot\,}{\,\cdot}$ means either $\inpd{\cdot\,}{\,\cdot}_{L^2}$ (AC--type)
or $\inpd{\cdot\,}{\,\cdot}_{H^{-1}}$ (CH--type).
Next, we  give  two convex splitting forms of $F(\phi)$
  as  discussed   in
\eqref{F_split_ln} and \eqref{F_split_nl}.
 The convex splitting form $F(\phi) = F_c^\l-F_e^\n$ can be taken as
\begin{equation}\label{F_split1}
 F_{c}^\l(\phi) = \int_0^1 \Big[ \frac{\kappa^2}{2}(\frac{\partial \phi}{\partial x})^2 + \phi^2 + \frac{1}{4} \Big]\,dx, \quad F_{e}^\n(\phi) = \int_0^1 -\frac{1}{4}\phi^4 + \frac{3}{2}\phi^2\, dx.
 \end{equation}
and the convex splitting form $F(\phi) = \titF_c^\n - \titF_e^\l$  is chosen as
\begin{equation}\label{titF_split}
 \titF_c^\n(\phi) = \int_0^1 \Big[ \frac{\kappa^2}{2}(\frac{\partial \phi}{\partial x})^2 + \frac{1}{4}\phi^4 + \frac{1}{4} \Big]\,dx, \quad \titF_e^\l(\phi) = \int_0^1 \frac{1}{2}\phi^2\, dx.
\end{equation}

\begin{remark}\label{rem:3}
The functionals in \eqref{titF_split}  are always convex
but  the  second functional $F_e^{\n}$ in \eqref{F_split1} are   convex
only in the region  $\phi\in [-1,1]$.
In general, when the global convexity is not available,
to obtain the locally contractive  $F_c$ and expansive $F_e$,
one usually has to introduce a sufficiently large positive constant $C$.
For instance the quadratic term  $-\phi^2$ in \eqref{eqn:F} can be written as $
 C\phi^2 -  ( 1/2+C )\phi^2$,
then the convex region of $F_e^{\n}$ is $ {\phi}^2 \leq (2C+1)/3$.
The form  \eqref{F_split1} corresponds to $C=1$.

The $L_2$ gradient flow  of the Ginzburg-Landau functional
ensures that the solution $\phi(t)$ always remains in  $[-1,1]$  by the maximum principle
if
the initial $\abs{\phi(0)} \leq 1$;
however, the $H^{-1}$ Cahn-Hillard flow may not always satisfy this condition.
It is not a trivial work to obtain such $L_\infty$ bounds
{\it a prior} and to derive  an  implementable optimal  choice of  $C$.
 Refer to   \cite{ElseyWirth2013}
   for  the  theoretic  investigation  on a convex splitting scheme for
a phase field crystal model.
In  our  numerical simulations,
we choose the minimal $C$ such that the local convexity holds for the initial
and it happens that we did not observe unstable phenomena from our empirical results of two applications
tested here. However, it should be noted that this choice has no theoretic foundation
and  is not guaranteed to work in any situation.
\end{remark}

\subsubsection{Saddle points in $L^2$ metric}\label{sec:saddle_in_L2}
The second-order variational operator $\Hess(\phi)$ of the energy functional $F$, evaluated at $\phi$,
is
 $ \Hess(\phi)\psi  =  \delta^2_\phi F \,  \psi= -\kappa^2 \Delta\psi   +  f''(\phi)\psi,
 ~~\forall\psi\in H^2([0,1]),
 $
 where $H^2$ is the standard Soblev space.
 The eigenvalue problem for this operator  is defined by
\begin{equation}\label{eig-AC}
\Hess(\phi)\psi = -\kappa^2 \Delta\psi   +  f''(\phi)\psi  =\lambda \psi,
\end{equation}
 subject to boundary conditions,
where $\lambda$ is the eigenvalue.
By the   result of the Rayleigh quotient,
the eigen-pair for the min-mode,   $\{\lambda_{1},\psi_1\}$,
solves the variational problem
\[ \min_{\psi\in H^1([0,1])}\mathcal{R}(\psi)
:= \frac{\inpd{\psi}{\Hess\psi}_{L^2}}{\|\psi\|^2_{L^2}} = \frac{  \int_0^1  {\kappa^2}  \abs{\nabla \psi}^2 + f''(\phi)\psi^2 \, dx }{\int_0^1 \abs{\psi}^2 dx }.
\]

\medskip

After the min-mode  is obtained, the subproblem of minimizing the auxiliary functional
\eqref{L-ACH} is then solved by evolving the gradient flow
 \begin{equation}\label{L2_metric}
\begin{split}
\frac{\partial \phi}{\partial t} &= -\frac{\delta L}{\delta \phi}(\phi)
= -\frac{\delta F}{\delta \phi}(\phi)
+2 ( v \otimes v ) {\frac{\delta F}{\delta \phi}(\hat\phi)},
\\
&= \kappa^2\Delta\phi - \phi^3 + \phi -2\vvT
( \kappa^2 \Delta\hat{\phi} -  \hat{\phi}^3 +\hat{\phi} ).
\end{split}
\end{equation}
$v=v^{(k+1)}$ is the min-mode of $\Hess(\phi^{(k)})$.
$\hat{\phi}$ is defined in \eqref{370}.

 \medskip

\emph{Convex Splitting Scheme.}
 We apply \eqref{F_split1} and \eqref{titF_split} to the convex splitting  of $L(\phi)$, then
$
 L(\phi)  =
\left[F_{c}^\l(\phi) + 2\,\titF_e^\l (\hat\phi)
\right]-
\left [F_{e}^\n(\phi) + 2\,\titF_c^\n(\hat\phi)
\right] =: L_c(\phi) - L_e(\phi),
 $
where
\begin{align*}
L_c(\phi) &= \int_0^1 \left[ \frac{\kappa^2}{2}({\phi_x})^2 + \phi^2 + \frac{1}{4} +
{\hat\phi}^2 \right] \, dx, \\
L_e(\phi) &= \int_0^1 \left[ -\frac{1}{4}\phi^4 + \frac{3}{2}\phi^2\ + \kappa^2 ({\hat\phi}_x)^2 + \frac{1}{2} {\hat\phi}^4 + \frac{1}{2} \right] \, dx.
\end{align*}
The first order variational derivatives of $L_c(\phi)$ and $L_e(\phi)$  are, respectively,
 \begin{align*}
 \delta_\phi L_c(\phi) &= -\kappa^2 \Delta\phi + 2\phi + 2(v\otimes v)\hat\phi, \\
 \delta_\phi L_e(\phi) &= -\phi^3 + 3\phi + 2(v\otimes v)\left(-\kappa^2\Delta{\hat\phi} + {\hat\phi}^3\right).
\end{align*}
Therefore, the  convex splitting scheme for (\ref{L2_metric}) is written as
 \begin{equation}\label{AC_convex_scheme1}
     \begin{split}
         \frac{\phi^{n+1}-\phi^n}{\Delta t}  =  & \left[ \kappa^2 \Delta\phi - 2\phi - 2(v\otimes v)\phi \right]^{n+1} \\
         & \!+\! \left[ -\phi^3 + 3\phi + 2(v\otimes v)\big( -\kappa^2\Delta\hat\phi + {\hat\phi}^3 \big) \right]^n,
     \end{split}
 \end{equation}
 after using $\langle v,v \rangle_{L^2} = 1$ and $(v\otimes v)\hat\phi = (v\otimes v)\phi$.

\medskip

\emph{Non-convex-splitting Scheme.} We need  a  traditional semi-implicit  scheme, which
is not derived from the  convex splitting idea, as the benchmark for comparison.
The idea is to take out the linearized part of the nonlinear term in $f'$ and
make it implicit together with  the Laplace operator.
Linearizing the cubic term $\phi^3$  in \eqref{L2_metric},
$
({\phi^{n+1}})^3 \approx ({\phi^n})^3 + 3 ({\phi^n})^2 (\phi^{n+1} - \phi^n)
$, we then have the following  scheme for \eqref{L2_metric}
\begin{align}\label{AC_nonconvex_scheme}
    \frac{\phi^{n+1} - \phi^n}{\Delta t} =
     &~ \kappa^2 \Delta\phi^{n+1}  - \left[ ({\phi^n})^3 + 3 ({\phi^n})^2 (\phi^{n+1} - \phi^n) \right] + \phi^{n+1} \nonumber\\
     & - 2 \left[ \langle v, \kappa^2 \Delta \hat\phi - {\hat\phi}^3 + \hat\phi \rangle v \right]^n.
\end{align}

\begin{figure}[h]
\begin{center}
\begin{subfigure}[b]{0.326\textwidth}
\includegraphics[width=\textwidth]{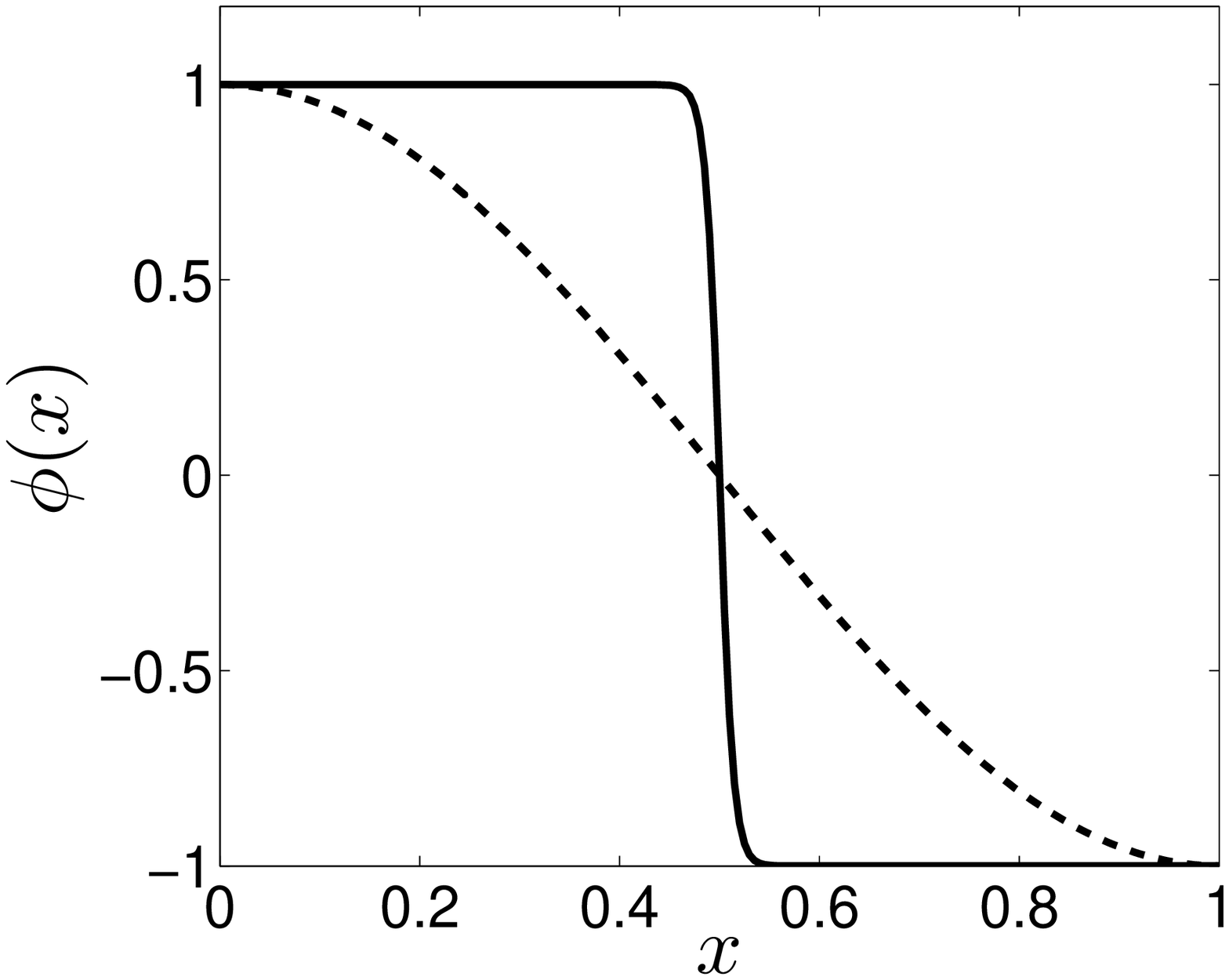}
\caption{}
\label{fig:AC:saddles}
\end{subfigure}
\hfill
\begin{subfigure}[b]{0.326\textwidth}
\includegraphics[width=\textwidth]{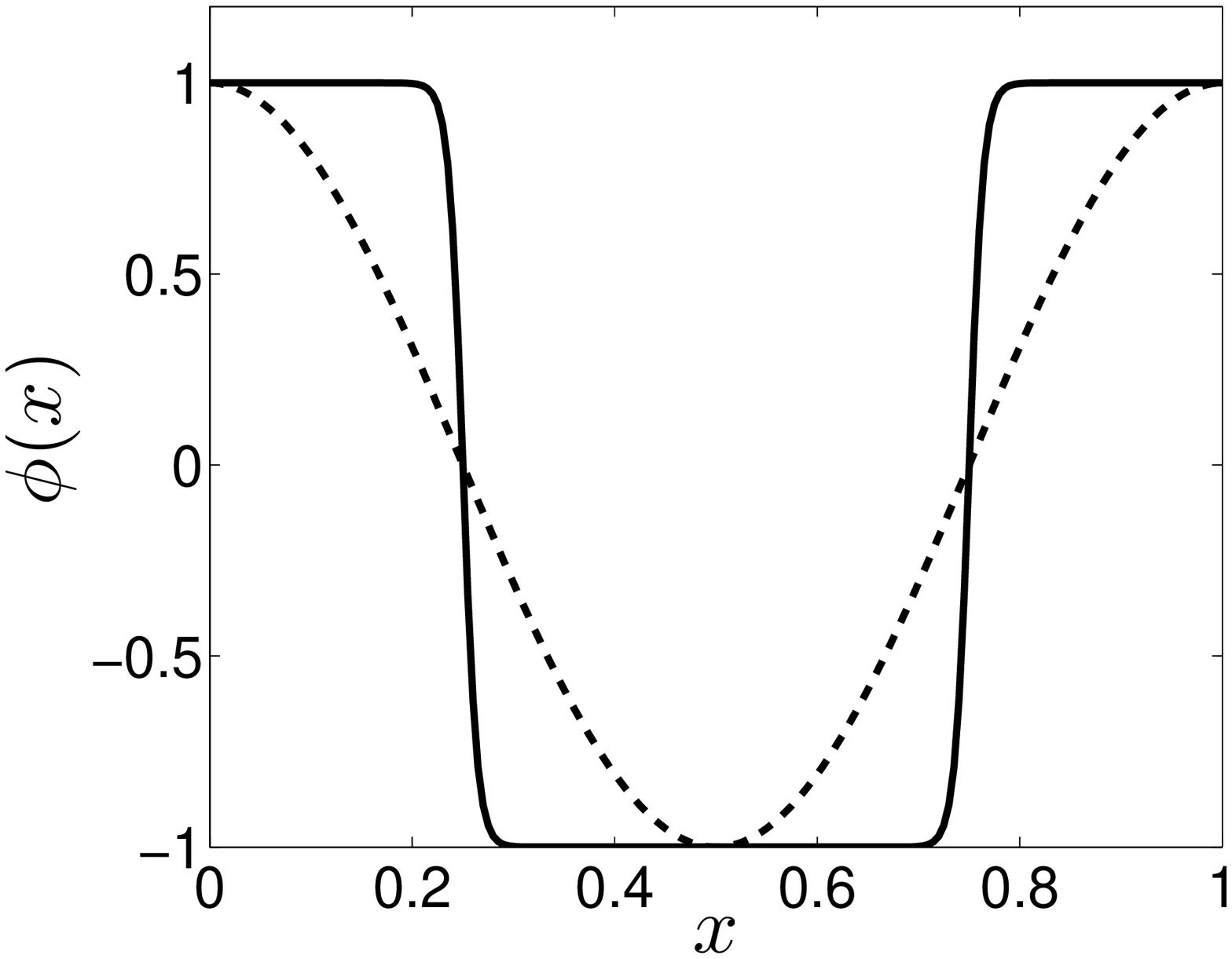}
\caption{}
\label{figA:AC_NBC:s1}
\end{subfigure}
\hfill
\begin{subfigure}[b]{0.326\textwidth}
\includegraphics[width=\textwidth]{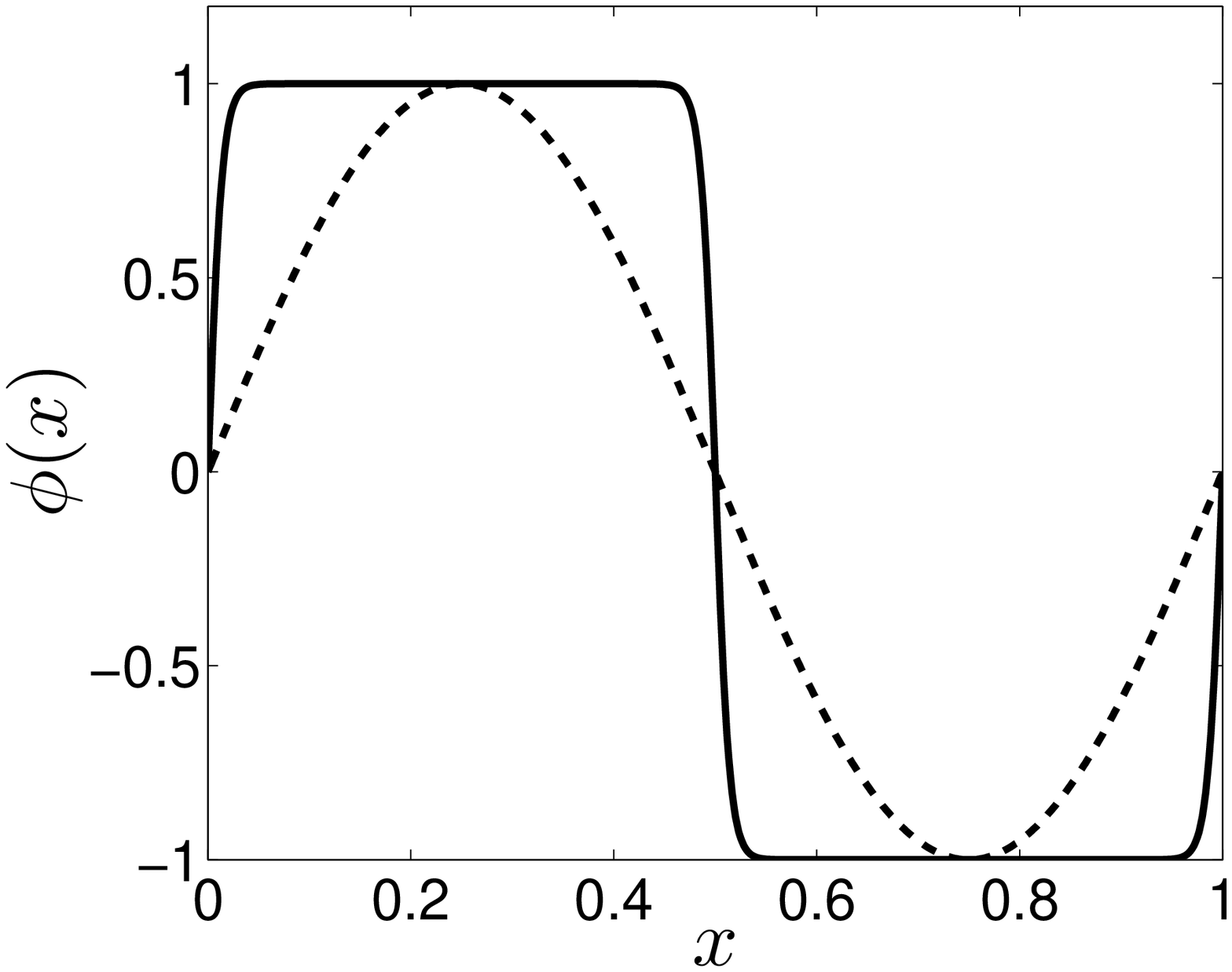}
\caption{}
\label{figB:AC_PBC:s2}
\end{subfigure}
\caption{ Profiles of some saddle points (solid lines) of $F(\phi)$ in $L^2$ metric  computed  from  various initial states
(dashed lines).  {\bf (a)}{\bf (b)}:   the  Neumann boundary condition;
{\bf (c)}:  the periodic boundary condition.
The free energy $F$ for these three plotted states from left to right are 0.0094, 0.0188 and 0.0188, respectively.
$\kappa=0.01$. }
\end{center}
\end{figure}

Set the parameter $\kappa=0.01$.
The finite difference method is used for spatial discretization with the mesh grid $\{x_i = i h , i=0, 1, 2,\ldots, N\}. ~h = 1/N.$ $N=200$.  A finer mesh with $N=1000$ is also used to verify all numerical results.
There are only  two locally stable states   in $L^2$ metric.
They are the two homogeneous constant states: $\phi_+(x)
\equiv 1$ and $\phi_-(x) \equiv -1$.
For the Neumann boundary condition, Figure \ref{fig:AC:saddles}
and Figure \ref{figA:AC_NBC:s1} show the transition states
calculated  from the convex splitting scheme \eqref{AC_convex_scheme1} with the initial conditions $\phi_0(x) = \cos{\pi x}$
and $\phi_0(x) = \cos 2\pi x$, respectively.
 For the periodic boundary condition,
 Figure \ref{figB:AC_PBC:s2} shows the    transition state
 obtained from the initial condition $\phi_0(x) = \sin 2\pi x$,
 which looks almost identical to the one in  Figure \ref{figA:AC_NBC:s1} after  a simple spatial translation.

\begin{figure}[htbp]
\begin{center}
\begin{subfigure}[b]{0.46\textwidth}
\includegraphics[width=\textwidth]{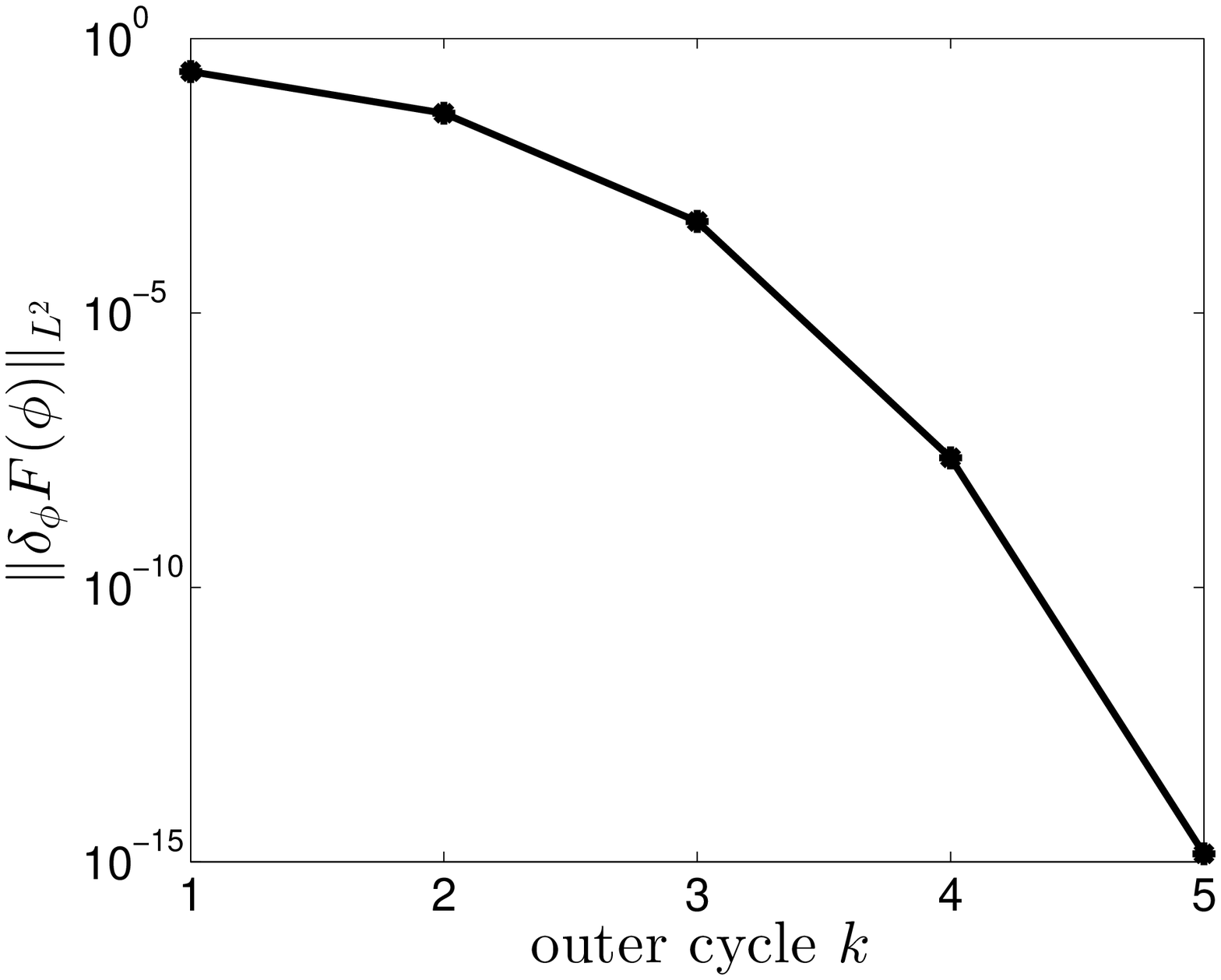}
\caption{Neumann boundary condition}
\label{fig:AC_NBC_rate}
\end{subfigure}
\hfill
\begin{subfigure}[b]{0.46\textwidth}
\includegraphics[width=\textwidth]{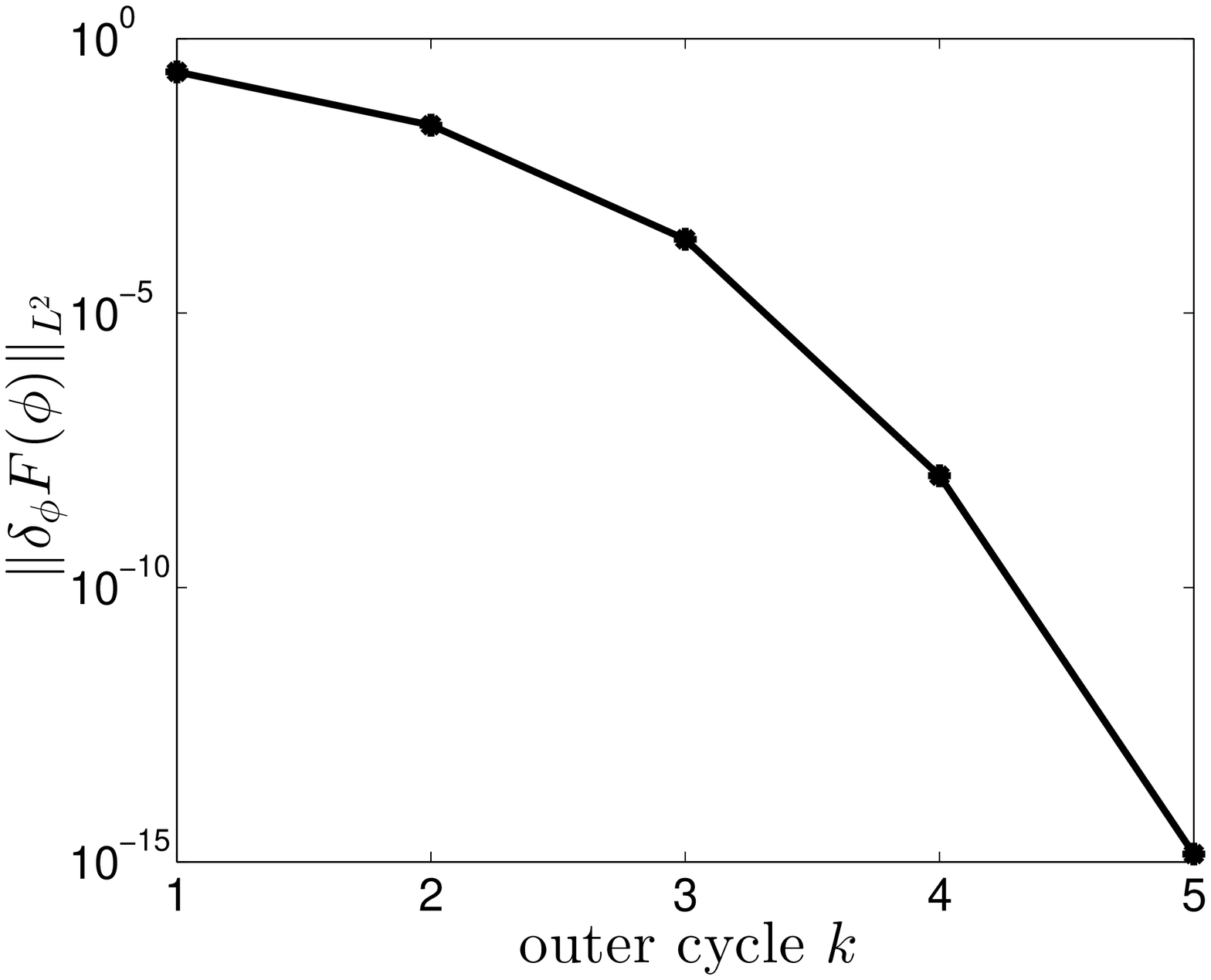}
\caption{Periodic boundary condition}
\label{figA:AC_PBC_conve_rate}
\end{subfigure}
\caption{The validation of the quadratic convergence rate of the IMF mapping $\Phi$
by plotting the decay of the error, measured by  the force $\|\delta_\phi F(\phi^{(k)})\|_{L^2}$ at each cycle $k$.
}
\label{fig:AC_PNBC_rate}
\end{center}
\end{figure}
Next we present  the numerical evidence on the performance of our numerical method. We first  validate  the quadratic convergence rate of the IMF for both boundary conditions.
 In   this validation, which only cares about the rate of the mapping $\phi\to \Phi(\phi)$,
 we actually solve each subproblem with an extremely  high precision.
 The numerical results are presented in Figure \ref{fig:AC_PNBC_rate}
 and are consistent with Theorem \ref{Th_IMF}.

\begin{table}[htbp]
  \begin{subtable}{\linewidth}
  \centering
  {
\begin{tabular}{|*{7}{r|}}
\hline
\multicolumn{7}{|c|}{The  required number of iterations } \\ \cline{1-7}
{\multirow{2}*{$\Delta t$}}{}  &  \multicolumn{2}{|c|}{$err = 10^{-4}$}  & \multicolumn{2}{|c|}{$err = 10^{-6}$} & \multicolumn{2}{|c|}{$err = 10^{-8}$} \\\cline{2-7}
{} & (\ref{AC_convex_scheme1})  &  (\ref{AC_nonconvex_scheme}) & (\ref{AC_convex_scheme1}) &  (\ref{AC_nonconvex_scheme}) & (\ref{AC_convex_scheme1}) &  (\ref{AC_nonconvex_scheme}) \\ \hline
0.01 & 777 & 763 & 1092 & 1074 & 1407 & 1386 \\ \hline
0.1 & 92 & 78 & 129  & 111 & 166 & 144 \\ \hline
  5.0 & 16 & $\infty$ & 23 &  $\infty$ & 30 & $\infty$ \\ \hline
  10 & 14 & $\infty$ & 21 & $\infty$ & 28 & $\infty$ \\ \hline
\end{tabular}
}
    \caption{Neumann boundary condition.
    The initial state is $\phi^{(0)}=\cos{\pi x}$. }
\end{subtable}
\vfill
  \begin{subtable}{\linewidth}
  \centering
  {
\begin{tabular}{|*{7}{r|}}
\hline
\multicolumn{7}{|c|}{The  required number of iterations } \\\cline{1-7}
{\multirow{2}*{$\Delta t$}}{}  &  \multicolumn{2}{|c|}{$err = 10^{-4}$}  & \multicolumn{2}{|c|}{$err = 10^{-6}$} & \multicolumn{2}{|c|}{$err = 10^{-8}$} \\\cline{2-7}
{} & (\ref{AC_convex_scheme1}) &  (\ref{AC_nonconvex_scheme}) & (\ref{AC_convex_scheme1}) &  (\ref{AC_nonconvex_scheme}) & (\ref{AC_convex_scheme1}) & (\ref{AC_nonconvex_scheme}) \\ \hline
0.01 & 728 & 722 & 1043 & 1034  & 1358 & 1347 \\ \hline
0.1 & 86 & 74 & 123 & 108 & 160 & 141 \\ \hline
5.0 & 14 & $\infty$ & 22 & $\infty$ & 29 & $\infty$ \\ \hline
  10 & 13 & $\infty$ & 20 & $\infty$ & 27 & $\infty$ \\ \hline
\end{tabular}
}
    \caption{Periodic  boundary condition.
    The initial state is $\phi^{(0)}=\sin{2\pi x}$. }
\end{subtable}

\vskip 2mm
\caption{The comparison of the CS scheme \eqref{AC_convex_scheme1} and nCS scheme \eqref{AC_nonconvex_scheme}  for the subproblem within the first  cycle $\phi^{(0)}\rightarrow \phi^{(1)}=\Phi(\phi^{(0)})$.
The integers shown in the table are
 the  required number of iterations   for the CS scheme (\ref{AC_convex_scheme1}) and the nCS (\ref{AC_nonconvex_scheme})
 to achieve the three choices of the prescribed  error tolerance $\|\delta_\phi L(\phi^n)\|_{L^2} \leq   10^{-4},   10^{-6}$ and $10^{-8}$.}\label{tab:AC_comparison}
\end{table}

Now  we compare the performance of
   the convex splitting   (``CS'') scheme
(\ref{AC_convex_scheme1}) against   the non-convex splitting  (``nCS'') scheme (\ref{AC_nonconvex_scheme}).
Firstly, we  examine  their performance for the subproblem, i.e., within  a fixed cycle
where only the inner iteration is  running.
Take the first  cycle for example, which is to solve $\min_{\phi} L(\phi; \phi^{(0)}, v^{(1)})$
with the initial $\phi(t=0)=\phi^{(0)}$.
We  measure the error by the gradient force $err:=\|\delta_\phi L(\phi^n)\|_{L^2}$
and calculate the iteration number  required
to attain the given error tolerance
for the CS scheme
\eqref{AC_convex_scheme1}  and the nCS scheme
\eqref{AC_nonconvex_scheme}.
The tolerances we tested are the following three values:
$err \leq 1.0 \times 10^{-4}, err\leq 1.0 \times 10^{-6}$ and $err\leq 1.0 \times 10^{-8}$. We can observe from Table \ref{tab:AC_comparison}
 for this simple subproblem  that (1) the CS scheme  obviously
 have much better stability than  the nCS scheme  when the   time step size is   large
 ($\infty$ in this table means that the numerical results   diverge);
 (2) For the small $\Delta t$ such that both schemes converge,   both schemes have the similar  required number of iterations (i.e., the time steps)
 to decrease $\delta_\phi L$
 to certain accuracy. This is expected since both solves the same steepest descent flow for $L$ with the same step size.

\begin{table}[htbp]
\begin{subtable}{\linewidth}
\centering
{
\begin{tabular}{|*{7}{r|}}
\hline
\multicolumn{7}{|c|}{ The   number of   cycles } \\\cline{1-7}
{\multirow{2}*{$\Delta t$}}
  {} &  \multicolumn{2}{|c|}{iter$ \#= 50 $}  & \multicolumn{2}{|c|}{iter$ \#= 80 $}  & \multicolumn{2}{|c|}{iter$\#= 100 $} \\\cline{2-7}
{} & \eqref{AC_convex_scheme1} & \eqref{AC_nonconvex_scheme} & \eqref{AC_convex_scheme1} & \eqref{AC_nonconvex_scheme} & \eqref{AC_convex_scheme1} & \eqref{AC_nonconvex_scheme} \\ \hline
0.01 & 29 & 28 & 18 & 18 & 15 & 14 \\ \hline
0.1 & 4 & 3 & 3 & 2 & 2 & 2 \\ \hline
5.0 & 1 & $\infty$ & 1 & $\infty$ & 1 & $\infty$ \\ \hline
\end{tabular}
}
 \caption{Neumann boundary condition.
    The initial state is $\phi^{(0)}=\cos{\pi x}$. }
    \end{subtable}
\vfill
\begin{subtable}{\linewidth}
\centering
{
\begin{tabular}{|*{7}{r|}}
\hline
\multicolumn{7}{|c|}{ The   number of   cycles } \\\cline{1-7}
{\multirow{2}*{$\Delta t$}}
  {} &  \multicolumn{2}{|c|}{iter$ \#= 50 $}  & \multicolumn{2}{|c|}{iter$ \#= 80 $}  & \multicolumn{2}{|c|}{iter$\#= 100 $} \\\cline{2-7}
{} & \eqref{AC_convex_scheme1} & \eqref{AC_nonconvex_scheme} & \eqref{AC_convex_scheme1} & \eqref{AC_nonconvex_scheme} & \eqref{AC_convex_scheme1} & \eqref{AC_nonconvex_scheme} \\ \hline
0.01 & 28 & 27 & 17 & 17 & 14 & 14  \\ \hline
0.1 & 4 & 3 & 3 & 2 & 2 & 2 \\ \hline
5.0 & 1 & $\infty$ & 1 & $\infty$ & 1 & $\infty$ \\ \hline
\end{tabular}
}
    \caption{Periodic  boundary condition.
    The initial state is $\phi^{(0)}=\sin{2\pi x}$. }
\end{subtable}
\caption{ The comparison of the  required number of outer cycles  for the CS  scheme \eqref{AC_convex_scheme1}
and the nCS scheme \eqref{AC_nonconvex_scheme}
to  attain  the given  error tolerance $\|\delta_\phi F(\phi^{(k)})\|_{L^2} \leq  10^{-8}$
  when   the inner iteration number is fixed as  50, 80 and 100, respectively,
  for different choices of the
   time step size.
 }\label{table:AC_whole_effect}
\end{table}

Secondly, we show  the comparison of
the  overall efficiency  in locating  the saddle point.
To be more transparent, we fix the number of iterations in each cycle
and count the required number of the outer cycles to reach some prescribed
tolerance for the error which is defined as $\norm{\delta_\phi F(\phi^{(k)})}_{L^2}$.
 The results are summarized in  Table \ref{table:AC_whole_effect}.
 The total iteration number, which is the indicator of the total computational cost,
 is therefore equal to the number of cycles multiplied by  the ``iter$ \#$'' specified in the corresponding columns of Table \ref{table:AC_whole_effect}.
The key conclusions from this table, for both Neumann and periodic boundary conditions, are the following:
(1) for the CS scheme, the  larger the time step size is, the smaller the total computational cost is; (2) for the nCS scheme, it is divergent when  a  large time step size is applied;
(3) for the very small time step size such as $0.01$, there is no much difference in
the computational cost of
the CS and nCS schemes. The last observation
is consistent with the known empirical experience of the convex splitting schemes:
when both the CS and nCS converge,
the decay rate of the objective function in the CS may not
be better.

\begin{figure}[htbp]
\begin{center}
\begin{subfigure}[b]{0.49\textwidth}
\includegraphics[width=\textwidth]{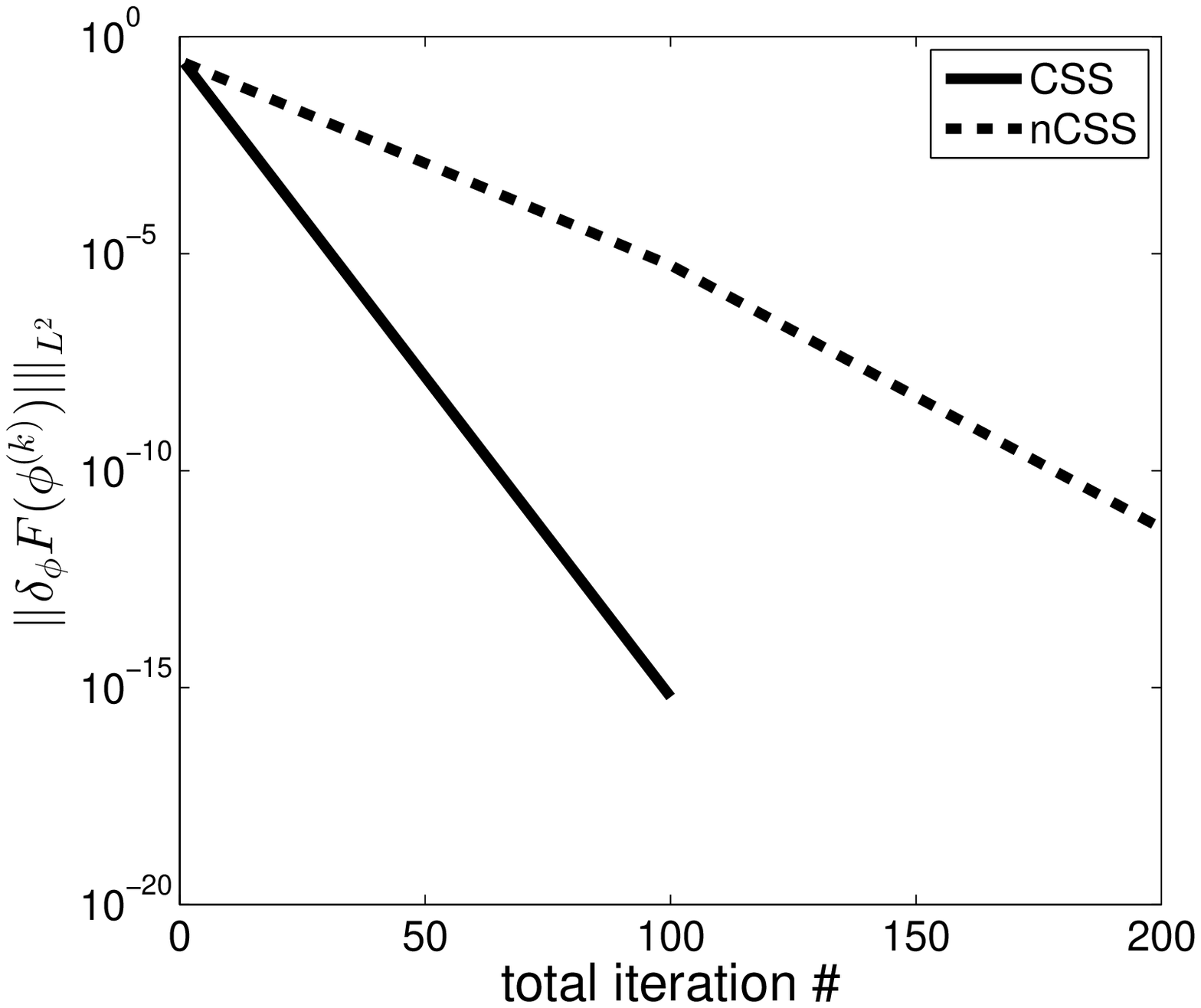}
\caption{Neumann boundary condition}
\label{FigA:error_AC_iterN_300}
\end{subfigure}
\begin{subfigure}[b]{0.49\textwidth}
\includegraphics[width=\textwidth]{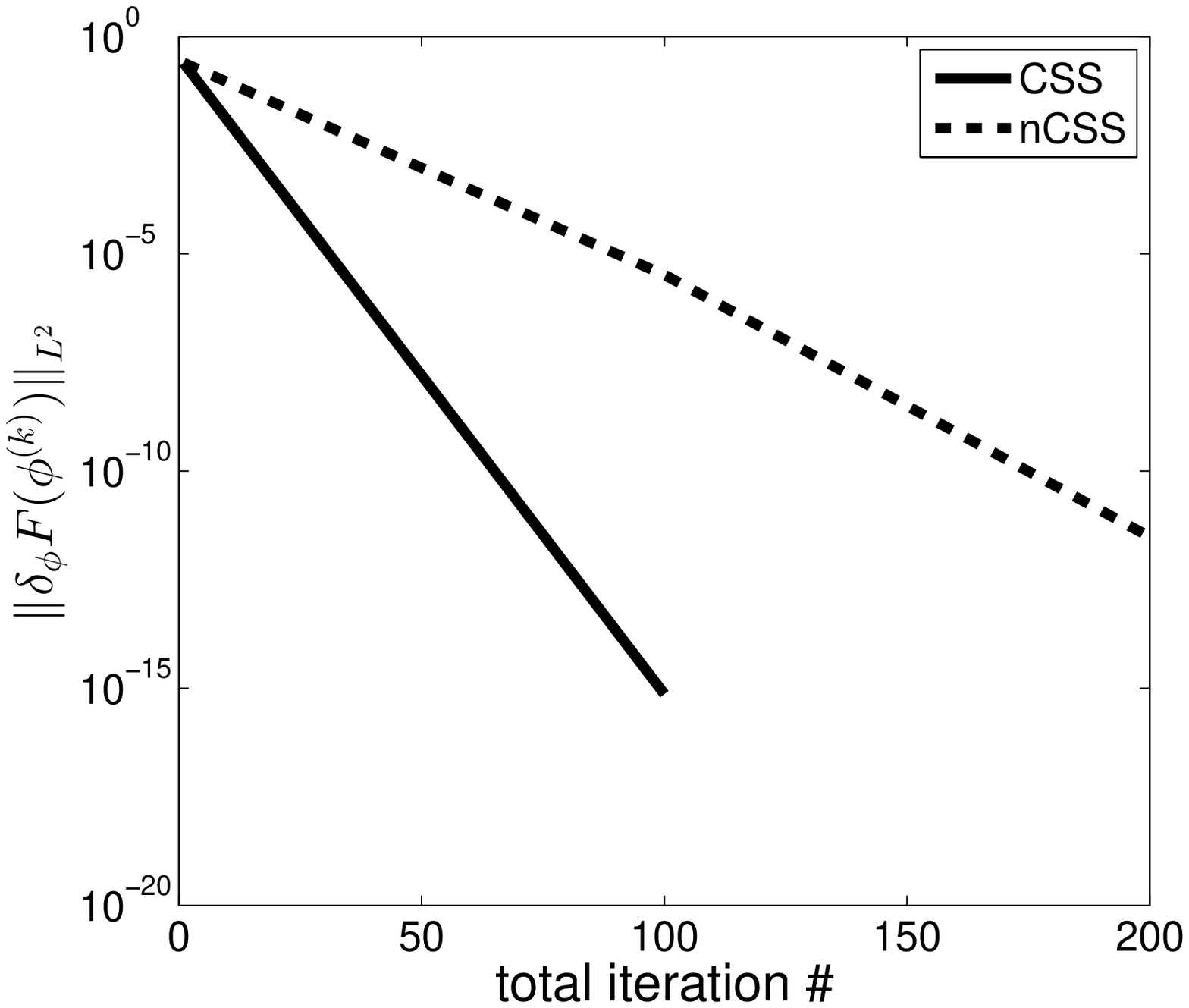}
\caption{Periodic boundary condition}
\label{figB:AC_PBC_compar}
\end{subfigure}
\caption{ The decay of the error measured by the gradient  $\|\delta_\phi F(\phi^{(k)})\|_{L^2}$ with the total iteration number
(i.e., the cost) for the CS scheme \eqref{AC_convex_scheme1} (solid line) and the nCS scheme \eqref{AC_nonconvex_scheme} (dashed line).  The inner iteration number in each cycle is fixed as 100 and  the  time step sizes are $\Delta t = 5.0 $ and $  0.1 $ for the CS and nCS schemes, respectively.}
\label{err_AC_compa}
\end{center}
\end{figure}

To better visualize the improvement of the CS scheme over  the nCS scheme,
we draw  the decay of the error $\|\delta_\phi F(\phi)\|_{L^2}$ with respect to the total iteration number in Figure \ref{err_AC_compa}.
 This plot   illustrates how much accuracy one can obtain (the vertical axis)
 with  the available total computational cost (the horizontal axis),
 in which the solid lines are from the CS scheme and the dashed lines are from the nCS scheme.
Note that the step sizes for the CS and nCS are different in order to obtain this
efficiency gain.

\subsubsection{Saddle points in $H^{-1}$ metric }
 In this part, we study  the transition state of $F(\phi)$ in $H^{-1}$ metric. Note that
the inner product and the norm in $H^{-1}$
metric can be transformed  to those in $L^2$ metric as follows
\begin{align}\label{-1to2}
 \norm{\phi}_{H^{-1}}^2 = \big\langle (-\Delta)^{-1}\phi,\phi \big\rangle_{L^2}, ~~
 \langle\phi,\psi\rangle_{H^{-1}} = \big\langle(-\Delta)^{-1}\phi,\psi \big\rangle_{L^2},
\end{align}
 where    $(-\Delta)^{-1}$,  a bounded positive self-adjoint linear operator,
 is the inverse of $-\Delta $ subject to certain boundary condition.
   Then the variational derivatives between
   the  $L^2$ metric and the $H^{-1}$ metric can be linked as follows:
$$ \delta_\phi F \Big|_{H^{-1}} = -\Delta \delta_\phi F, \quad
\wt{\Hess}:=  \delta^2_\phi F \Big|_{H^{-1}} = -\Delta \delta^2_\phi F, $$
 where $\delta_\phi F$ and $\delta^2_\phi F$ represent respectively the first and the second order variational derivatives of $F(\phi)$ in $L^2$ metric.

The gradient flow   $\partial_t \phi= -( -\Delta \delta_\phi F)=-\kappa^2 \Delta^2 \phi + \Delta (\phi^3-\phi)$ is the CH equation \eqref{CH_eq}. It is known that
the solution $\phi(x,t)$ of  the CH equation \eqref{CH_eq}
satisfies  $\int_0^1 \phi(x,t)dx \equiv \int_0^1 \phi(x,0)dx, \forall t>0 $.
In our problem of finding the saddle point in the $H^{-1}$ metric,
we choose  a fixed mass $m=0.6$   beforehand and we are interested in the saddle points satisfying  $\int_0^1 \phi (x) dx = m $.
We also require any   eigenvectors or perturbations  to belong to the
subspace
$\set{\psi: \int_0^1 \psi (x) dx =0}$.
The eigenvalue problem of  $\wt{\Hess}$  is

\begin{equation}\label{eig-CH}
\begin{cases}
\wt{\Hess}(\phi)\psi  =-\Delta (-\kappa^2 \Delta\psi   +  f''(\phi)\psi) =\lambda \psi,
\\
\qquad    \int_0^1 \psi(x) \,dx = 0,
\end{cases}
\end{equation}
subject to the Neumann or periodic boundary condition.
For nonzero eigenvalue, the eigenvector $\psi$ automatically satisfies
$\int_0^1 \psi dx =0$; For zero eigenvalue, the condition $\int_0^1 \psi dx =0$ needs to be imposed additionally.
If we introduce the projection ${\mathbf{P}}u:= u - \int_0^1 u(x)\, dx $, then the eigen-problem \eqref{eig-CH}
is equivalent to   ${\mathbf{P}}\wt{\Hess}{\mathbf{P}}u = \lambda  u$
for any $u$ without the zero-mass constraint.  The min-mode is then equal to $u$ if
the   eigenvalue is nonzero and    equal to $\mathbf{P}u$ if
the   eigenvalue is zero.

For the periodic boundary condition,
there exists one degenerate ($\lambda=0$) direction $\psi_0$ for the Hessian $\wt{\Hess}(\phi)$,
which is the spatial derivative of $\phi$, i.e.,
\[ \psi_0(x)=\phi_x(x) .\]
To see this, we just need verify \eqref{eig-CH} by using the periodic boundary condition:
$\int_0^1 \phi_x(x) d x = \phi(1) -\phi(0)=0$ and
\[
\int_0^1 -\Delta (-\kappa^2 \Delta\phi_x   +  f''(\phi)\phi_x)  dx
=\int_0^1 \kappa^2 \partial_x^5\phi(x) - \partial_x^2 (  \partial_x (f'(\phi(x)))) ~~dx=0.\]

 The Rayleigh quotient with respect to $H^{-1}$ metric  is
\[
\wt{\mathcal{R}} (\psi)= \frac{\inpd{\psi}{ \wt{\Hess}\psi}_{H^{-1}} }{\norm{\psi}^2_{H^{-1}}}
=\frac{\int_0^1 \kappa^2 \abs{\nabla \psi}^2+ f''(\phi)\psi^2 \, dx}{\int_0^1 \psi \Delta^{-1} \psi dx},
\]
and thus  the min-mode is  the minimizer
of the problem
\begin{equation}\label{R-CH}
  \argmin_{\psi} \set{\wt{\mathcal{R}} (\psi):  {\int_0^1 \psi \, dx =0, ~~\norm{\psi}_{H^{-1}}=1} }.
\end{equation}

For the IMF in  this $H^{-1}$ case,
the subproblem of minimizing the  auxiliary functional $L$  for a given $\phi^{(k)}$ at cycle $k$ is
\begin{equation}\label{eqn:CHk}
 \phi^{(k+1)}=\argmin_{\int_0^1 \phi(x) dx=m}  L(\phi; \phi^{(k)}, v^{(k+1)}),
 \end{equation}
where $m=\int_0^1 \phi^{(k)} dx $.
This means that each IMF cycle $\phi^{(k)}\to \phi^{(k+1)}$
should also conserve the mass.
The expression of $L(\phi)$ is defined  as \eqref{L-ACH}, with the modification of  $\hat\phi$ as follows
  \begin{equation} \label{1371}
  \begin{split} \hat{\phi} & := \phi^{(k)} + \inpd{v}{\phi-\phi^{(k)}}_{H^{-1}} v
   = \phi^{(k)} + \inpd{ -\Delta ^{-1} v}{\phi-\phi^{(k)}}_{L^2} v\\
  & = \phi^{(k)} + \inpd{w}{\phi-\phi^{(k)}}_{L^2} v,
    \end{split}
\end{equation}
{  where }~~~
\[ w: =-\Delta^{-1} v
\]
is the unique solution satisfying the equation  $-\Delta w =v$  and
 $\int_0^1 w\, dx = 0$.
Then, $
\delta_\phi L(\phi) = \delta_\phi F(\phi) - 2 w \inpd{\delta_\phi F(\hat\phi)}{v}_{L^2},
$
and the $H^{-1}$ gradient flow of $L$   is
\begin{equation}\label{H_1metric}
 \begin{split}
  \frac{\partial\phi}{\partial t} &= \Delta\frac{\delta L}{\delta\phi}(\phi)
= \Delta \left(\frac{\delta F}{\delta\phi}(\phi) \right) + 2  v   \inpd{  \frac{\delta F}{\delta\phi}(\hat{\phi} ) }{v}_{L^2},
\\
 & = -\kappa^2 \Delta^2\phi + \Delta(\phi^3 - \phi) + 2v \inpd{v}{-\kappa^2\Delta\hat{\phi}+\hat{\phi}^3-\hat\phi}_{L^2},
 \end{split}
 \end{equation}
where $v=v^{(k+1)} $ refers to the min-mode of  \eqref{eig-CH} at $\phi^{(k)}$
(normalized under $H^{-1}$ metric, i.e.,  $\norm{v}_{H^{-1}}=1$).
Note that the scalar  in the  second line of \eqref{H_1metric}
is    the $L^2$ inner product due to the cancelation of
$\Delta$ and $\Delta^{-1}$ in the calculation, but the variable
$\hat{\phi}$   involves the $H^{-1}$ metric via $w$.
This suggests one more computational cost
of solving $w=-\Delta^{-1}v$ during the IMF subproblem
than the GAD scheme where $\hat{\phi}$ is actually
just $\phi^{(k)}$ (\cite{IMA2015}).
Our result  \eqref{H_1metric}
is consistent with the equation (3.4)
in \cite{MMSLZZ2013}  which was written  for the GAD in the finite dimension, i.e.,
 $\hat{\phi}=\phi^{(k)}$.

\begin{remark}
\label{rem:4}
We show that the  flow   \eqref{H_1metric}  conserves the
initial mass $\int_0^1 \phi(x) dx $,  sharing   exactly  the  same property  as  the CH equation \eqref{CH_eq}. So the constraint in \eqref{eqn:CHk}
holds automatically.
This result immediately implies  that the IMF mapping
$\phi^{(k)}\to\phi^{(k+1)} $ does not change the mass at each cycle $k$. To prove our conclusion,
after  integrating the two sides of \eqref{H_1metric} and using the boundary conditions
(either Neumann or periodic), one remains to  show the following condition for the eigenvector $v$:
$\int_0^1 v(x) dx =0.$
This is exactly the condition that the min-mode $v$  satisfies in \eqref{eig-CH}.
\end{remark}

\medskip

\emph{Convex splitting scheme.}
We now  test the  convex splitting form of  \eqref{F_split1},\eqref{titF_split}, which
corresponds to
$ L(\phi) = L_c(\phi) - L_e(\phi),$
where
$ L_c(\phi) = F_{c}^l(\phi) + 2\tilde F_e^l(\hat\phi)$,
$ L_e(\phi) = F_{e}^n(\phi) + 2\tilde F_c^n(\hat\phi)$,
and $\hat{\phi}$ is defined in \eqref{1371}.
Since $
 \delta_\phi L_c(\phi) = -\kappa^2 \Delta\phi + 2\phi + 2 \inpd{v}{\hat\phi}_{L^2} w, $
 and $
 \delta_\phi L_e(\phi) = -\phi^3 + 3\phi + 2\inpd{-\kappa^2\Delta\hat\phi + \hat\phi^3}{v}_{L^2}w
$,
then the convex splitting scheme for \eqref{H_1metric} is for $n\geq 0$
\begin{equation}\label{CH_convex_scheme1}
    \begin{split}
        \frac{\phi^{n+1}-\phi^n}{\Delta t} = & \left[ -\kappa^2 \Delta^2\phi^{n+1} + 2\Delta\phi^{n+1} - 2\inpd{w}{\phi^{n+1}}_{L^2}
        \inpd{v}{v}_{L^2}  v \right] \\
        & - 2\inpd{v}{\phi^{(k)}}_{L^2}v+  2 \inpd{w}{\phi^{(k)}}_{L^2} \inpd{v}{v}_{L^2}v \\
        & + \left[ \Delta (\phi ^{n})^3 - 3\Delta\phi ^{n} + 2 \inpd{v}{-\kappa^2 \Delta\hat{\phi}^n + (\hat{\phi}^n)^3}_{L^2} v \right],
    \end{split}
\end{equation}
where $v=v^{(k+1)}$, $w=-\Delta^{-1}v^{(k+1)}$ and $\hat{\phi}^n $ is  from \eqref{1371} by letting $\phi=\phi^n$.
On the right hand side of  this CS scheme \eqref{CH_convex_scheme1}, the first line
is linear in the unknown $\phi^{n+1}$, the second line is independent of $n$,
the third line is the nonlinear term in $\phi^n$.

\emph{Non-convex-splitting scheme.}
The same linearization trick   used in \eqref{AC_nonconvex_scheme} for  the previous AC type problem
is   applied here again to construct the following  non-convex splitting scheme for CH type equation for comparison:
\begin{equation}\label{CH_nonconvex_scheme}
    \begin{split}
        \frac{\phi^{n+1}-\phi^n}{\Delta t} = & -\kappa^2 \Delta^2\phi^{n+1} - \Delta\phi^{n+1} +  \Delta \left[   ({\phi^n})^3 +
        3 ({\phi^n})^2 (\phi^{n+1} - \phi^n) \right] \\
        & -  2 \left[ \langle v, \kappa^2 \Delta \hat\phi - {\hat\phi}^3 + \hat\phi \rangle v \right]^n.
    \end{split}
\end{equation}

 The parameter $\kappa=0.04$.
  We set the mesh grid $\{x_i = i h , i=0, 1, 2,\ldots, N\}. ~h = 1/N.$ $N=200$.
 The function $\phi(x)$ is represented by
  $\bm{\phi} = (\phi_0, \phi_1, \phi_2, \cdots, \phi_N)^T$ for the Neumann boundary condition
  and   $\bm{\phi} = (\phi_0, \phi_1, \phi_2, \cdots, \phi_{N-1})^T$ for the periodic boundary condition; $\phi_i \approx \phi(x_i)$.
 The matrix form of the Hessian $\wt{\Hess}$  at the state $\bm{\phi}$
 is
 $  \wt{\Hess}(\bm{\phi})=A\nabla^2  F_h(\bm{\phi}),$
 where $\nabla^2 F_h(\bm{\phi}) = \frac{\partial^2 F_n}{\partial\phi_i \partial\phi_j}$.  $A$ and $F_h(\bm\phi)$ denote  the discretized forms
  of the operator $-\Delta$ and the potential energy $F(\phi)$, respectively.
 The min-mode $v$ of the Hessian matrix $\wt{\Hess}$ can be calculated
 according to  the Rayleigh quotient \eqref{R-CH}.

\begin{figure}[htbp]
\begin{center}
\begin{subfigure}[b]{0.46\textwidth}
\includegraphics[width=\textwidth]{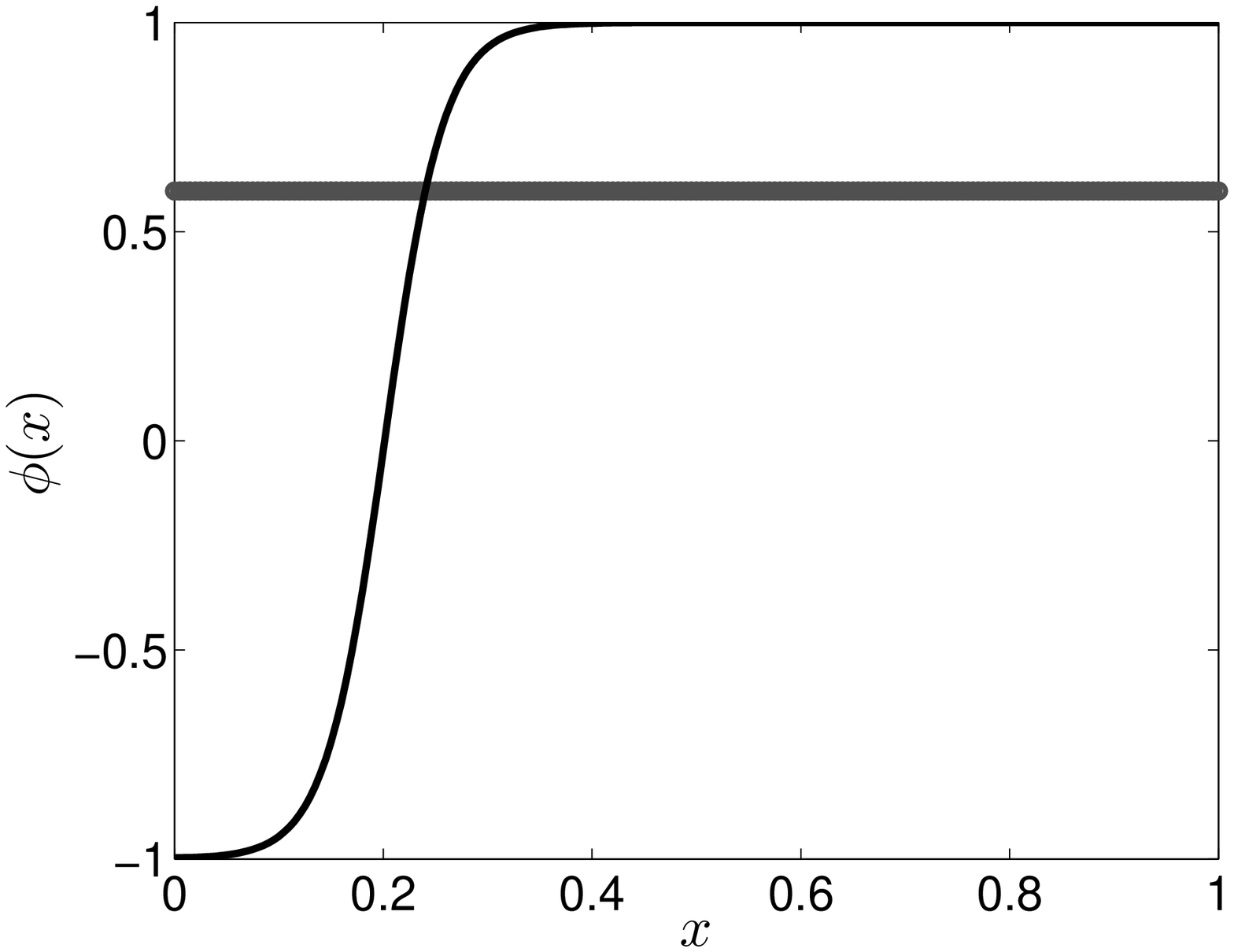}
\caption{Two locally stable states}
\label{figA:CH_NBC:metastable states}
\end{subfigure}
\hfill
\begin{subfigure}[b]{0.46\textwidth}
\includegraphics[width=\textwidth]{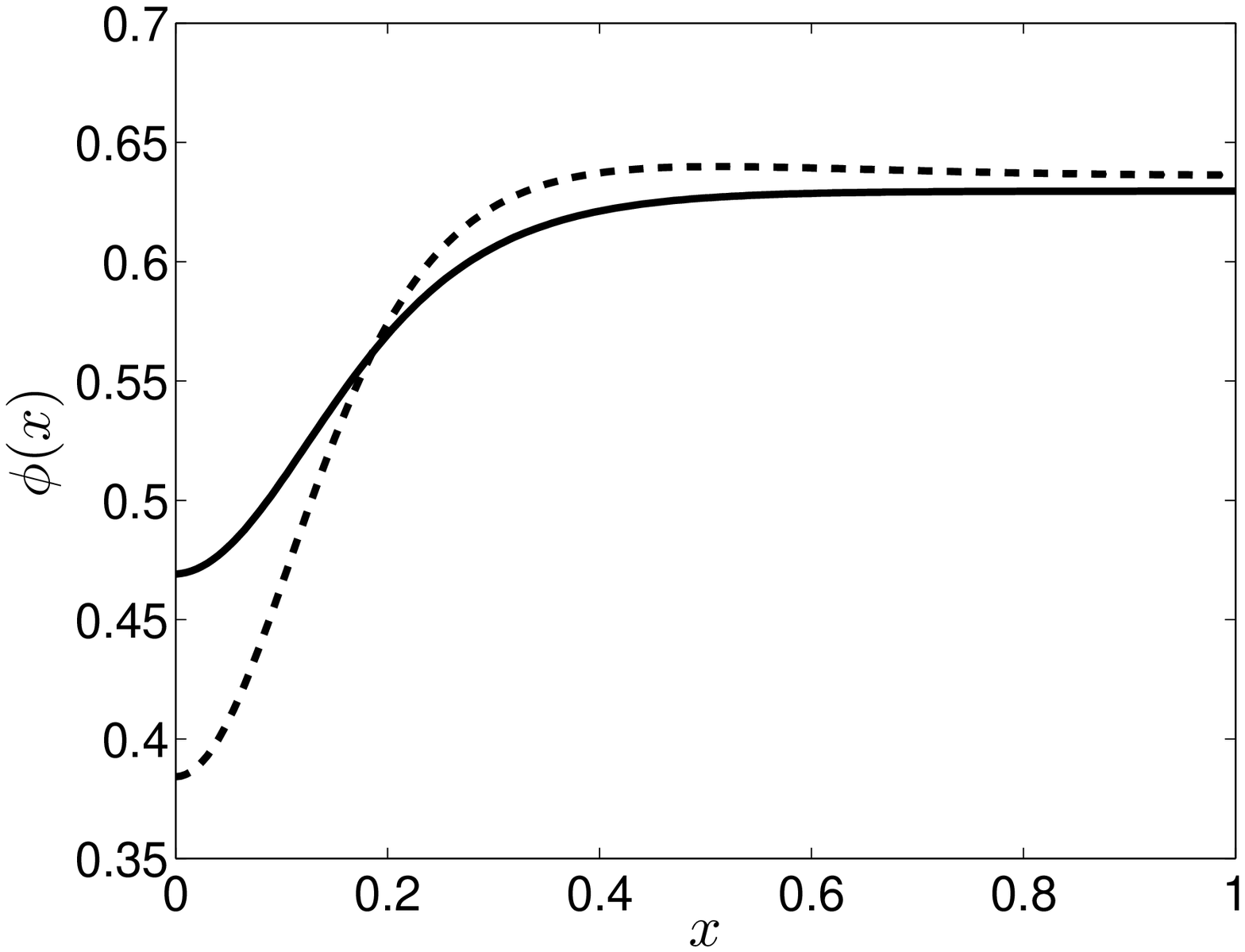}
\caption{Transition state}
\label{figB:CH_NBC:transition state}
\end{subfigure}
\caption{(Neumann boundary condition.   $H^{-1}$ metric.) ~{\bf (a)}: the  two
 stable stationary states of $F(\phi)$  with the  mass $\int_0^1 \phi\, dx =0.6$.
$F=0.10240$ for the trivial constant   state (the thick line) and $F = 0.03772$ for the transition layer state (the thin line).
{\bf (b)}: one of the
   transition states (solid line) with the free energy $0.10241$
   whose  first $3$ eigenvalues are $\lambda=-3.41, 3.91$ and 18.14,
   calculated from the initial condition (dashed line) whose first $3$ eigenvalues are $\lambda = -10.97, 3.45$ and 17.48.
Note that the   vertical axes in   subfigures are    in   different  scales.}
\label{fig:CH:fixed points_NBC}
\end{center}
\end{figure}

\begin{figure}[htbp]
\begin{center}
\begin{subfigure}[b]{0.45\textwidth}
\includegraphics[width=\textwidth]{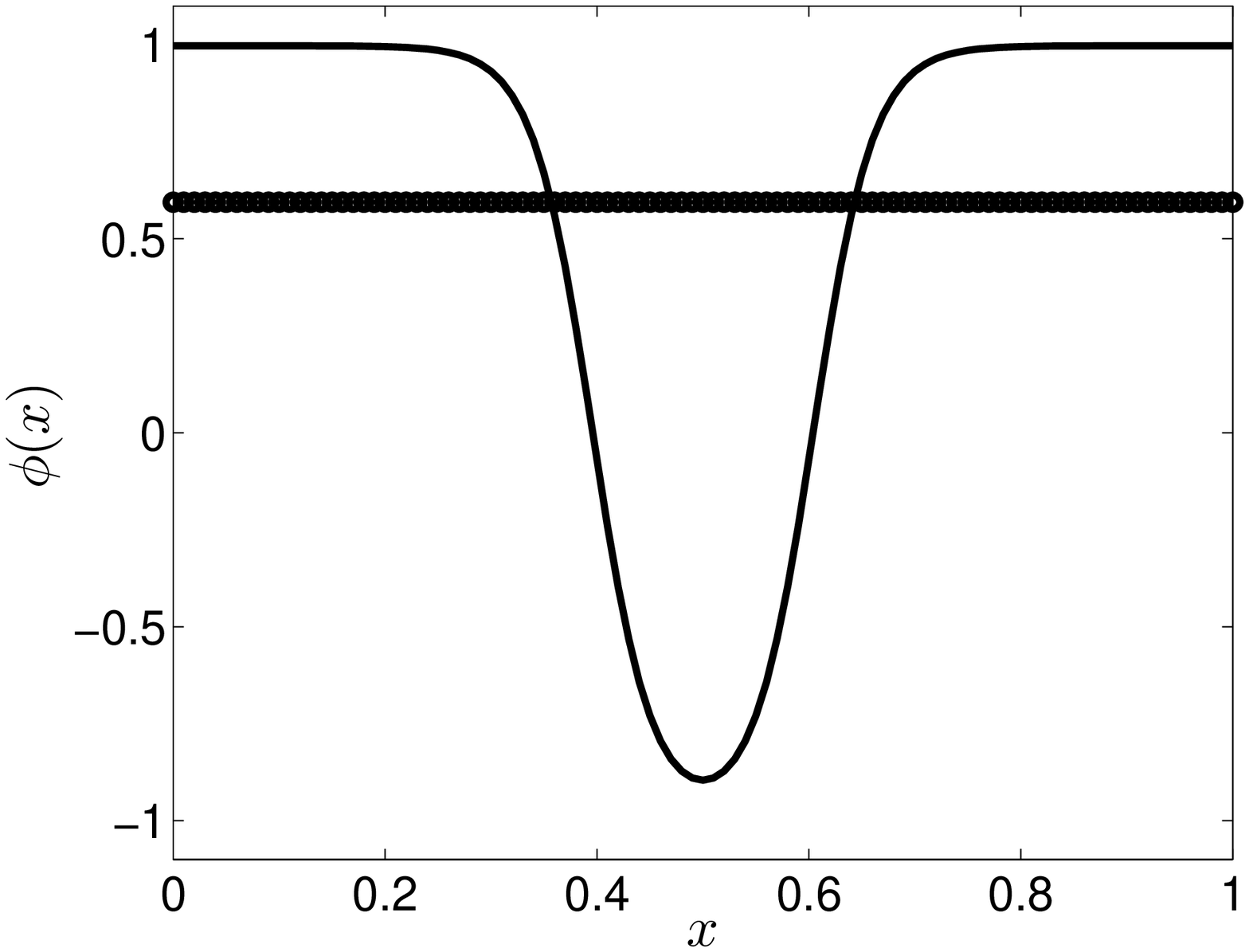}
\caption{Two locally stable states}
\label{figA:CH:metastable states}
\end{subfigure}
\hfill
\begin{subfigure}[b]{0.45\textwidth}
\includegraphics[width=\textwidth]{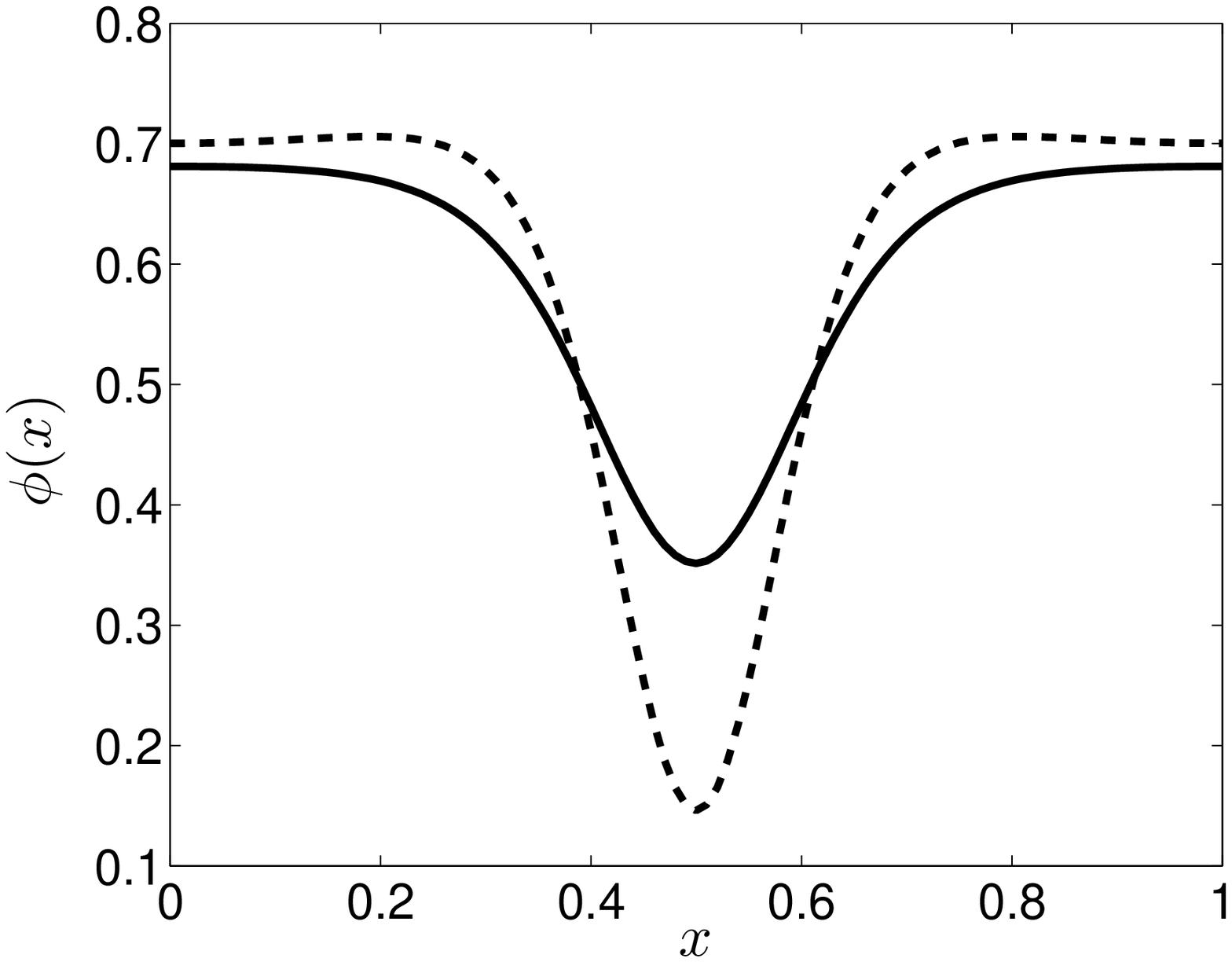}
\caption{Transition state}
\label{figB:CH:transition state}
\end{subfigure}
\caption{(Periodic boundary condition.  $H^{-1}$ metric.) ~ {\bf (a)}:  the  two stable stationary states of $F(\phi)$ in
  $H^{-1}$ metric with the given   mass $\int_0^1 \phi\, dx =0.6$.
$F=0.10240$ for the constant   state and $F=0.07510$ for the other
non-constant state.
{\bf (b)}:  one of the transition states (solid line) with the free energy $0.10285$ whose first three smallest (numerical) eigenvalues are $-12.75, -1.86\times 10^{-7}$ (corresponding to
theoretical zero eigenvalue) and $46.15$.
The initial state is shown as the dashed line in which the first 3 smallest eigenvalues are $\lambda = -27.13, -6.81$ and 44.98.
}
\label{fig:CH:fixed points}
\end{center}
\end{figure}
Figure \ref{fig:CH:fixed points_NBC}
and Figure \ref{fig:CH:fixed points}
show some of the  stationary states  of $F$ in    $H^{-1}$ metric
for the Neumann and periodic boundary conditions, respectively.
The stationary  points identified by us agree   with the result    in \cite{CH1D2012}.
Next we  compare the performance of the CS scheme \eqref{CH_convex_scheme1} and the nCS shceme \eqref{CH_nonconvex_scheme}.
 We also start from their performance for the subproblem in the first cycle as before. Table \ref{tab:CH_NBC_comparison} compares
 the number of inner iterations to  attain  the  three  error tolerances: $ \|\Delta \delta_\phi L \|_{H^{-1}} \leq   10^{-4},  10^{-5}$ or $  10^{-6}$.   Basically we have the same observations as for the $L_2$-metric case. The CS scheme allows the step size as large as $0.1$ but the nCS is not able in this example.
This table also shows that the nCS scheme, if it converges,   requires a less number of iterations to reach the tolerance, in particular for the periodic boundary condition.  This is possible because a scheme with a better stability does not guarantee a faster convergence rate toward the minimum.

\begin{table}[htbp]
\begin{subtable}{\linewidth}
\centering
{
\begin{tabular}{|*{7}{r|}}
\hline
\multicolumn{7}{|c|}{The required number of iterations} \\\cline{1-7}
{\multirow{2}*{$\Delta t$}}{}  &  \multicolumn{2}{|c|}{$err = 10^{-4}$}  & \multicolumn{2}{|c|}{$err = 10^{-5}$} & \multicolumn{2}{|c|}{$err = 10^{-6}$} \\\cline{2-7}
{} & (\ref{CH_convex_scheme1}) &  (\ref{CH_nonconvex_scheme}) & (\ref{CH_convex_scheme1}) &  (\ref{CH_nonconvex_scheme}) & (\ref{CH_convex_scheme1}) & (\ref{CH_nonconvex_scheme}) \\ \hline
 $10^{-3}$ & 2956 & 2647 & 3332 & 2993 & 3709 & 3339 \\ \hline
$10^{-2}$ & 550 & 273 & 617 & 308 & 685 & 343 \\
\hline
$10^{-1}$ & 283 & $\infty$ & 321 & $\infty$ & 360 & $\infty$ \\
\hline
$1.0$ & 14 & $\infty$ & 262 & $\infty$ & 298 & $\infty$ \\
\hline
\end{tabular}
}
\caption{ Neumann boundary condition.}
\end{subtable}
\vfill
\begin{subtable}{\linewidth}
\centering{
\begin{tabular}{|*{7}{r|}}
\hline
\multicolumn{7}{|c|}{The required number of iterations} \\\cline{1-7}
{\multirow{2}*{$\Delta t$}} {}  &  \multicolumn{2}{|c|}{$err = 10^{-4}$}  & \multicolumn{2}{|c|}{$err = 10^{-5}$} & \multicolumn{2}{|c|}{$err = 10^{-6}$} \\\cline{2-7}
{} & (\ref{CH_convex_scheme1}) & (\ref{CH_nonconvex_scheme}) & (\ref{CH_convex_scheme1}) & (\ref{CH_nonconvex_scheme}) & (\ref{CH_convex_scheme1}) & (\ref{CH_nonconvex_scheme}) \\ \hline
$10^{-3}$ & 1022 & 820 & 1163 & 939 & 1305 &  1057 \\ \hline
$10^{-2}$ & 281 & 77 & 316 & 88 & 351 & 99 \\ \hline
$10^{-1}$  & 208 &$ \infty $& 233 &$ \infty $& 257 &$ \infty$ \\ \hline
 \end{tabular}
}
\caption{Periodic boundary condition. }
\end{subtable}
\caption{The comparison of the CS  scheme \eqref{CH_convex_scheme1} and nCS scheme  \eqref{CH_nonconvex_scheme} for the subproblem $\phi^{(0)}\rightarrow \phi^{(1)}=\Phi(\phi^{(0)})$.
The integers   in the table are
 the  required number of iterations
 to achieve the three   prescribed    tolerances $\|\Delta\delta_\phi L(\phi^n) \|_{H^{-1}} \leq   10^{-4},
 10^{-5}$ and $  10^{-6}$.  }\label{tab:CH_NBC_comparison}
\end{table}

\begin{figure}[htbp]
\begin{center}
\begin{subfigure}[b]{0.50\textwidth}
\includegraphics[width=\textwidth]{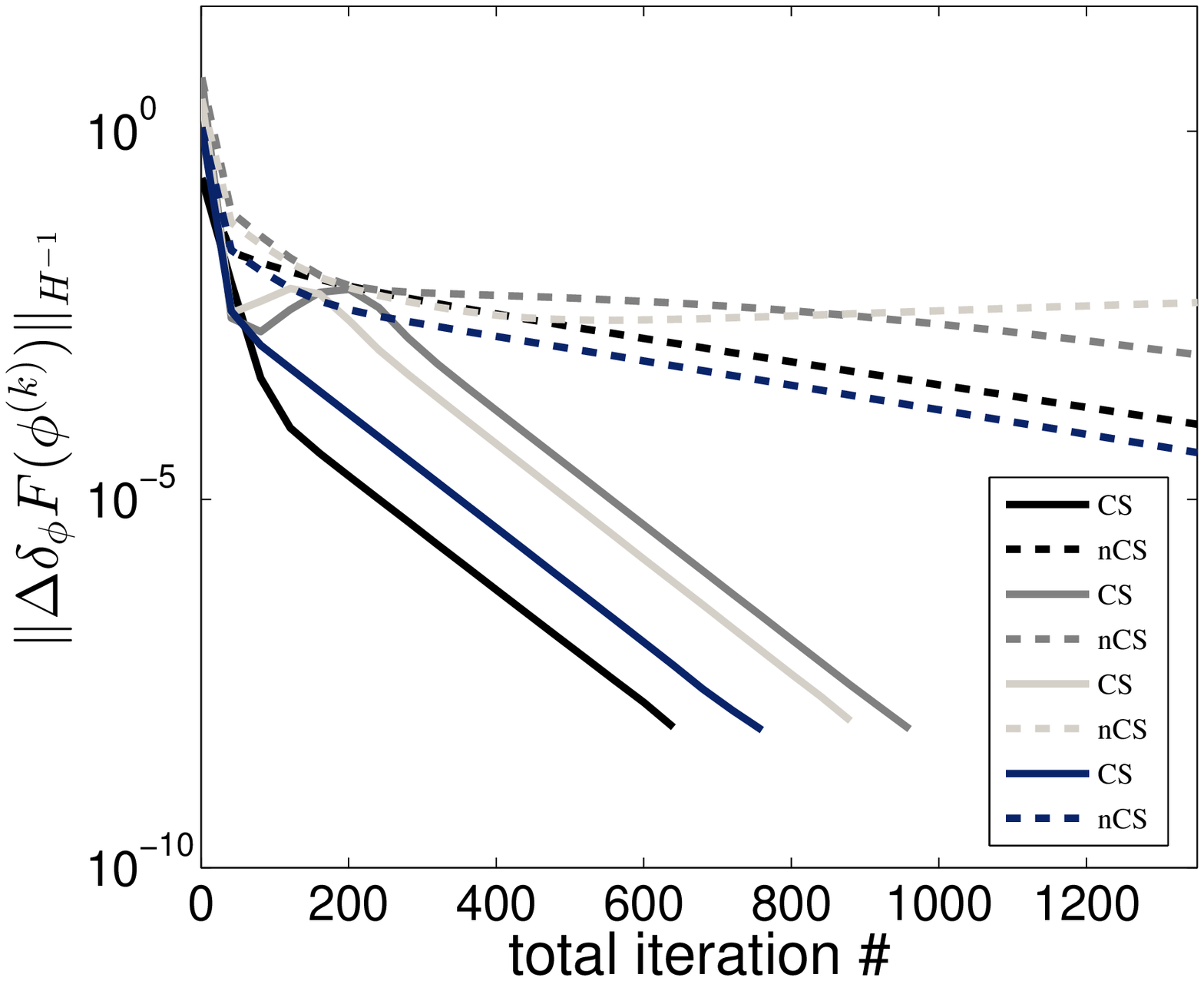}
\caption{Neumann boundary condition}
\label{figA:CH_NBC_compar}
\end{subfigure}
\begin{subfigure}[b]{0.49\textwidth}
\includegraphics[width=\textwidth]{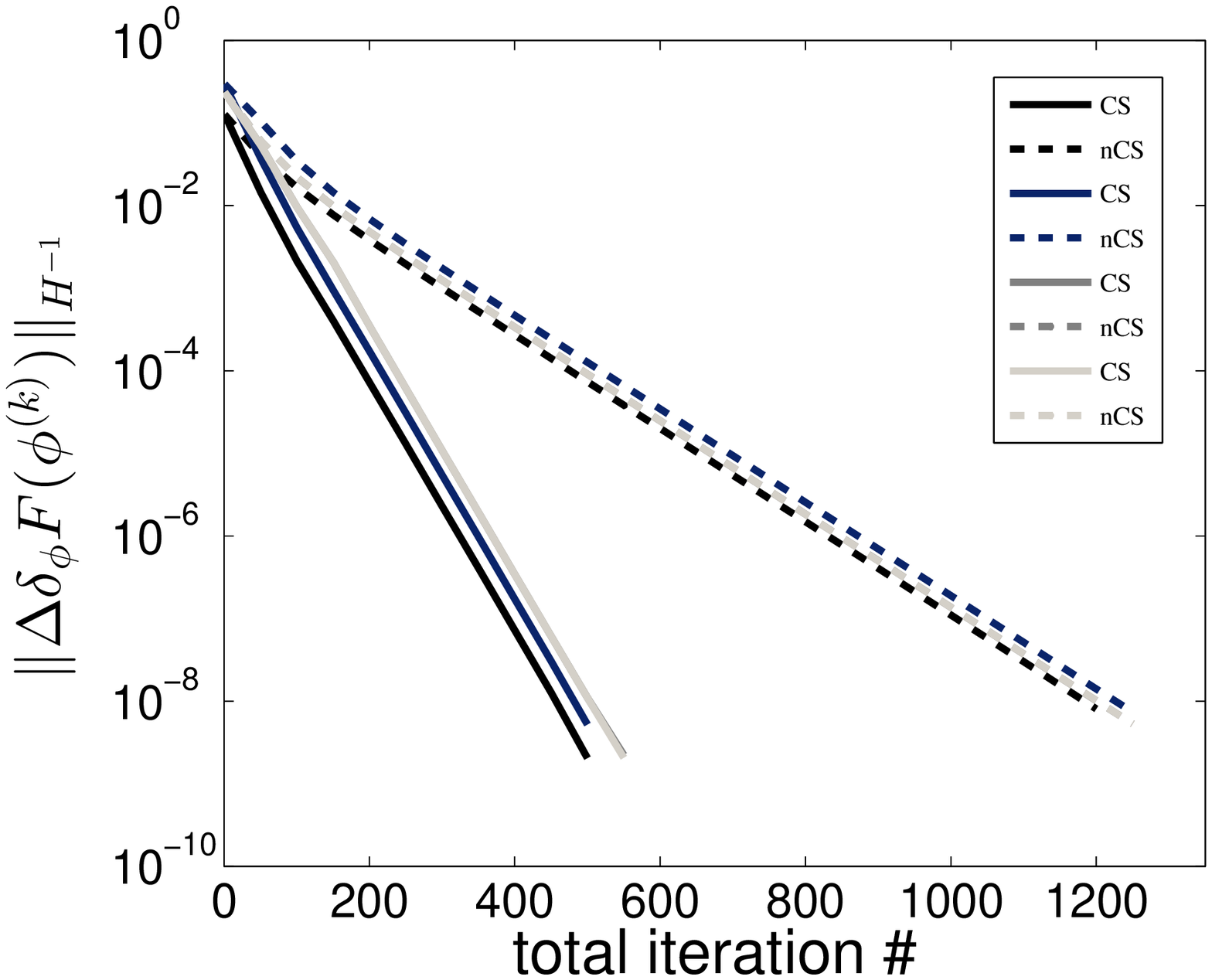}
\caption{Periodic boundary condition}
\label{figB:CH_PBC_compar}
\end{subfigure}
\caption{ The evolution of the errors measured by  $\|\Delta\delta_\phi F(\phi)\|_{H^{-1}}$
w.r.t.  the total iteration number for the CS scheme (\ref{CH_convex_scheme1}) (solid lines) and the nCS scheme (\ref{CH_nonconvex_scheme}) (dashed lines). Different small perturbations
around the initial condition used in Table  \ref{tab:CH_NBC_whole_effect2}
are added as the new initial conditions here to produce multiple lines. For the Neumann boundary condition in ({\bf a}),  the inner iteration number is  40 and  the  time step sizes are $\Delta t = 10^{-1} $ and $  10^{-3}$ for the CS   and nCS schemes, respectively. For the periodic boundary condition in ({\bf b}),
 the inner iteration number is 50 and the  time step sizes are the same as  in ({\bf a}).
}\label{CH_rate_compar}
\end{center}
\end{figure}

Table \ref{tab:CH_NBC_whole_effect2} shows
the same quantities of   the overall performance of the convex splitting scheme
in this $H^{-1}$ metric
as
 measured before in  Table  \ref{table:AC_whole_effect} in the $L_2$ metric.
 The conclusions we draw  from this table are qualitatively the same as
 from Table  \ref{table:AC_whole_effect}.
In the end, based on the experiments corresponding to   Table \ref{tab:CH_NBC_whole_effect2}, we plot the evolution of the error measured by the force $\|\Delta\delta_\phi F(\phi^{(k)}) \|_{H^{-1}}$
against  the total iteration number   for the CS scheme and the nCS scheme by using their own optimal time step sizes respectively.
Refer to  Figure \ref{CH_rate_compar}.
Here for the CH-type problem, in order to illustrate that our results are robust with respect to the initial guess,
we   added the random perturbation to generate multiple initial states so that multiple lines
are plotted for   different  initial guesses.

\begin{table}[htbp]
\begin{subtable}{\linewidth}
\centering{
\begin{tabular}{|*{7}{r|}}
\hline
\multicolumn{7}{|c|}{The   number of   cycles} \\\cline{1-7}
{\multirow{2}*{$\Delta t$}}
  {}  & \multicolumn{2}{|c|}{iter$ \#= 40 $}  & \multicolumn{2}{|c|}{iter$\#=50 $} & \multicolumn{2}{|c|}{iter$\#=60 $} \\\cline{2-7}
{} & (\ref{CH_convex_scheme1}) &  (\ref{CH_nonconvex_scheme}) & (\ref{CH_convex_scheme1}) &  (\ref{CH_nonconvex_scheme}) & (\ref{CH_convex_scheme1}) & (\ref{CH_nonconvex_scheme}) \\ \hline
$10^{-3}$ & 118 & 102 & 95 & 82 & 79 & 69 \\ \hline
 $10^{-2}$ & 28 & $\infty$ & 22 & $\infty$ & 16 & $\infty$ \\ \hline
  $10^{-1}$ & 16 & $\infty$ & 17 & $\infty$ & 16 & $\infty$ \\ \hline
\end{tabular}}
\caption{ Neumann boundary condition. }
\end{subtable}
\begin{subtable}{\linewidth}
\centering{
\begin{tabular}{|*{7}{r|}}
\hline
\multicolumn{7}{|c|}{The number of   cycles} \\ \cline{1-7}
{\multirow{2}*{$\Delta t$}}
  {} & \multicolumn{2}{|c|}{iter$ \#= 50 $}  & \multicolumn{2}{|c|}{iter$\#=80 $} & \multicolumn{2}{|c|}{iter$\#=100 $} \\\cline{2-7}
{} & (\ref{CH_convex_scheme1}) &  (\ref{CH_nonconvex_scheme}) & (\ref{CH_convex_scheme1}) &  (\ref{CH_nonconvex_scheme}) & (\ref{CH_convex_scheme1}) & (\ref{CH_nonconvex_scheme}) \\ \hline
$10^{-3}$ & 34 & 25 & 22 & 16 & 18 & 13 \\ \hline
 $10^{-2}$ & 13 & $\infty$ & 7 & $\infty$ & 8 & $\infty$ \\ \hline
 $10^{-1}$ & 11 & $\infty$ & 8 & $\infty$  & 7 & $\infty$ \\ \hline
\end{tabular}
}
\caption{Periodic boundary condition. }
\end{subtable}
 \caption{ The comparison of the  number of  outer cycles required for the CS scheme \eqref{CH_convex_scheme1} and the nCS
 scheme \eqref{CH_nonconvex_scheme} to attain the given error tolerance $\|\Delta\delta_\phi F(\phi^{(k)})\|_{H^{-1}} \leq  10^{-8}$, when    the inner iteration number is fixed for the different choices of the time step size $\Delta t=   10^{-3},  10^{-2}$ and $ 10^{-1}$.
 The corresponding initial states are specified in Figure \ref{figB:CH_NBC:transition state} and
 Figure \ref{figB:CH:transition state}, respectively.}
\label{tab:CH_NBC_whole_effect2}
\end{table}

{
\subsection{2D example: Landau-Brazovskii free energy}\label{sec4_2}
In this section, we study the nucleation problem of phase transition in diblock
copolymers (\cite{Tiejun2D,RobertA}).
The model is described by the Landau-Brazovskii energy functional
 of the order parameter $\phi$
\begin{equation}
F(\phi) = \int_\Omega   \frac{\xi^2}{2}[(\Delta+1)\phi(\textbf{r})]^2 + \Phi(\phi) ~d\,\textbf{r},
\end{equation}
where $\Phi(\phi) = \frac{\tau}{2}\phi^2 - \frac{\gamma}{3!}\phi^3 + \frac{1}{4!}\phi^4$.
The parameters are $\tau=-0.15, \xi=1.0, \gamma = 0.25$.
We compute the transition state of this $F$
 in $H^{-1}$ metric.

We consider  the two dimensional domain $\Omega = [0,L_x]\times [0,L_y]$ and the  periodic boundary condition.
 $\phi$ satisfies the mass conservation $\int_\Omega \phi(\textbf{r}) d\textbf{r} = 0$.
The   eigenvalue problem in the IMF for this case then reads
\begin{equation*}
\begin{cases}
\wt{\Hess}(\phi)\psi  = - \Delta \Big[ \xi^2 (\Delta+1)^2  +  \Phi''(\phi) \Big] \psi =\lambda \psi,
\\
\qquad    \int_\Omega \psi(\textbf{r}) \,d\textbf{r} = 0.
\end{cases}
\end{equation*}
 The min-mode $v$ is the eigenvector of $\wt{\Hess}$ corresponding to the smallest eigenvalue.
The  gradient flow associated with the minimization subproblem
for  $L(\phi; \phi^{(k)}, v^{(k+1)})
=F(\phi)- 2 F(\hat\phi)$
 is
 \begin{equation}\label{2D_L_eq}
 \frac{\partial\phi}{\partial t} = \Delta\frac{\delta L}{\delta\phi} = \xi^2\Delta (\Delta+1)^2\phi + \Delta \Phi_1(\phi) + 2 \inpd{v} { \xi^2 (\Delta+1)^2 \hat\phi + \Phi_1(\hat\phi)}_{L^2} v,
 \end{equation}
 where $\Phi_1(\phi) = \Phi'(\phi) = \tau\phi - \gamma\phi^2/2!+\phi^3/3! $ and $\hat\phi$ is defined the same as \eqref{1371}.

\medskip
 \emph{Convex-splitting scheme.}
We  give two convex splitting forms of $F$ as follows
 \begin{align}
 F_{c}^\l(\phi) &= \int_\Omega \frac{\xi^2}{2}\big[ (\Delta+1)\phi \big]^2 + \frac{\tau+22.75}{2} \phi^2 \,d\r, \label{2DF_split1}\\
 F_{e}^\n(\phi) &= \int_\Omega \big[ -\frac{1}{4!}\phi^4 + \frac{\gamma}{3!}\phi^3 + \frac{22.75}{2}\phi^2 \big]\, d\r,\label{2DF_split2}
 \end{align}
 and
 \begin{align}
 \titF_c^\n(\phi) &= \int_\Omega \frac{\xi^2}{2}|\Delta\phi|^2 + (\frac{1}{4!}\phi^4 - \frac{\gamma}{3!}\phi^3 + \frac{\xi^2+\tau}{2}\phi^2)\,d\r, \label{2DtitF_split1}\\
 \titF_e^\l(\phi) &= \int_\Omega \xi^2|\nabla\phi|^2 \, d\r.\label{2DtitF_split2}
 \end{align}

By \eqref{2DF_split1}, \eqref{2DF_split2} and \eqref{2DtitF_split2}, \eqref{2DtitF_split2},
 and  the convex splitting form of $L(\phi)=L_c(\phi)-L_e(\phi)$ with
$L_c(\phi) = F_{c}^\l(\phi) + 2 \titF_e^\l(\hat\phi)$ and $ L_e(\phi) = F_{e}^\n(\phi) + 2 \titF_c^\n(\hat\phi)$,
we have the convex splitting scheme for \eqref{2D_L_eq} as follows
\begin{align*}
\frac{\phi^{n+1}-\phi^n}{\Delta t} = \Delta [ \delta_\phi L_c(\phi) ]^{n+1} - \Delta [\delta_\phi L_e(\phi)]^n,
\end{align*}
i.e.,
\begin{align}\label{CSS_2D}
\frac{\phi^{n+1}-\phi^n}{\Delta t} = & \Big[ \xi^2\Delta(\Delta+1)^2 \phi + (\tau+22.75)\Delta\phi \Big]^{n+1} +  4\xi^2 \Big[ \inpd{\Delta v}{v} \inpd{w}{\phi} v \Big]^{n+1} \nonumber\\
& + 4\xi^2 \Big[\!\inpd{\Delta\phi^{(k)}}{v} \!-\! \inpd{w}{\phi^{(k)}} \inpd{\Delta v}{v} \! \Big] v \!-\! \Delta \Big[ \!-\frac{1}{3!} \phi^3 \!+\! \frac{\gamma}{2}\phi^2 + 22.75\phi \Big]^n \nonumber \\
& + 2 \inpd{v}{\xi^2\Delta^2\hat\phi + (\xi^2+\tau)\hat\phi - \frac{\gamma}{2}\hat\phi^2 + \frac{1}{3!}\hat\phi^3 }v \Big]^n,
\end{align}
where $\inpd{\cdot}{\cdot}$ means $\inpd{\cdot}{\cdot}_{L^2}$.
\begin{remark}
\label{rem:5}
In \eqref{2DF_split2},
we choose the constant  $C=22.75$ in Remark \ref{rem:3}, then   the convex region for \eqref{2DF_split2} is $[-6.5,6.5]$.
One can see later that the local minimum and the saddle point lie between $-1.5$ and $1.5$.
But our initial guess spans the interval  $[-6.5,6.5]$, so we use a large $C=22.75$.
Fortunately, we find that $\phi$ always locates in  $[-6.5,6.5]$ by tracking the numerical value of   $\phi$.
 Note that \eqref{2DF_split1}, \eqref{2DtitF_split1} and \eqref{2DtitF_split2} are globally convex.
 \end{remark}

\emph{Non-convex-splitting scheme.}
This  scheme is constructed by the same idea as before by applying the Taylor expansion of the nonlinear term $\Phi_1(\phi^{n+1})$ around the solution at the old time step $\phi^n$,
$$
\Phi_1(\phi^{n+1}) \approx \Phi_1(\phi^n) + \Phi^{ \prime}_1(\phi^n)(\phi^{n+1}-\phi^n).
$$
After simplification, we   get the non-convex splitting scheme as follows:
\begin{equation}\label{nCSS_2D}
    \begin{split}
        \frac{\phi^{n+1}-\phi^n}{\Delta t} =  &~ \xi^2 \Delta \left(\Delta+1\right)^2  \phi^{n+1} +  \Delta
         \left(\tau-\gamma\phi^n+\frac{1}{2}[{\phi^n}]^2\right) \phi^{n+1} \\
        & +   \Delta \left(\frac{\gamma}{2}[{\phi^n}]^2 - \frac{1}{3}{[\phi^n]}^3\right) + 2 \left[ \langle v, \xi^2 (\Delta+1)^2 \hat\phi + \Phi_1(\hat\phi) \rangle v \right]^n.
    \end{split}
\end{equation}

In the numerical simulation, we choose the domain $\Omega = [0,\frac{16}{\sqrt{3}}\pi]\times [0,8\pi]$ and set the mesh gird  $\{ x_i = ih_x, i=0,1,2,\ldots,N_x. \} $ and $\{ y_j = jh_y, j=0,1,2,\ldots,N_y\}$.  $h_x = 1/N_x, h_y = 1/N_y, N_x = N_y = 100.$
We first consider the gradient flow
$
\partial_t \phi = \Delta\delta_\phi F = \Delta[\xi^2(\Delta+1)^2\phi+\Phi_1(\phi)]
$
to find the steady states of $F(\phi)$.
 The initial conditions are chosen as  equation (2.21) in \cite{RobertA}.
For this steepest descent flow of $F$,
we tested  the convex splitting method based on    \eqref{2DF_split1}   \eqref{2DF_split2}
(obtained by setting $v=0$ in \eqref{CSS_2D}),
 and  the semi-implicit scheme based on  the Taylor expansion for the nonlinear term
 (this scheme is obtained  by directly setting $v=0$   in \eqref{nCSS_2D} ).
 The numerical  results show that all schemes  are stable  with   the time step sizes as large as $50$.

 Figure \ref{fig:2D_metastable_states}(a) and Figure \ref{fig:2D_metastable_states}(c) are   the two metastable states
 of the lamellar phase and the cylindrical phase respectively.
Figure \ref{fig:2D_metastable_states}(b) is the transition state obtained by the IMF.
The   initial guess for the IMF in our test
has the gradient   $\|\Delta\delta_\phi F(\phi^0)\|_{H^{-1}} = 10^{-2}$  and the minimal  eigenvalue    $-2.0\times 10^{-12}$.

\begin{figure}[htbp]
\begin{center}
\begin{subfigure}[b]{0.32\textwidth}
\includegraphics[width=\textwidth]{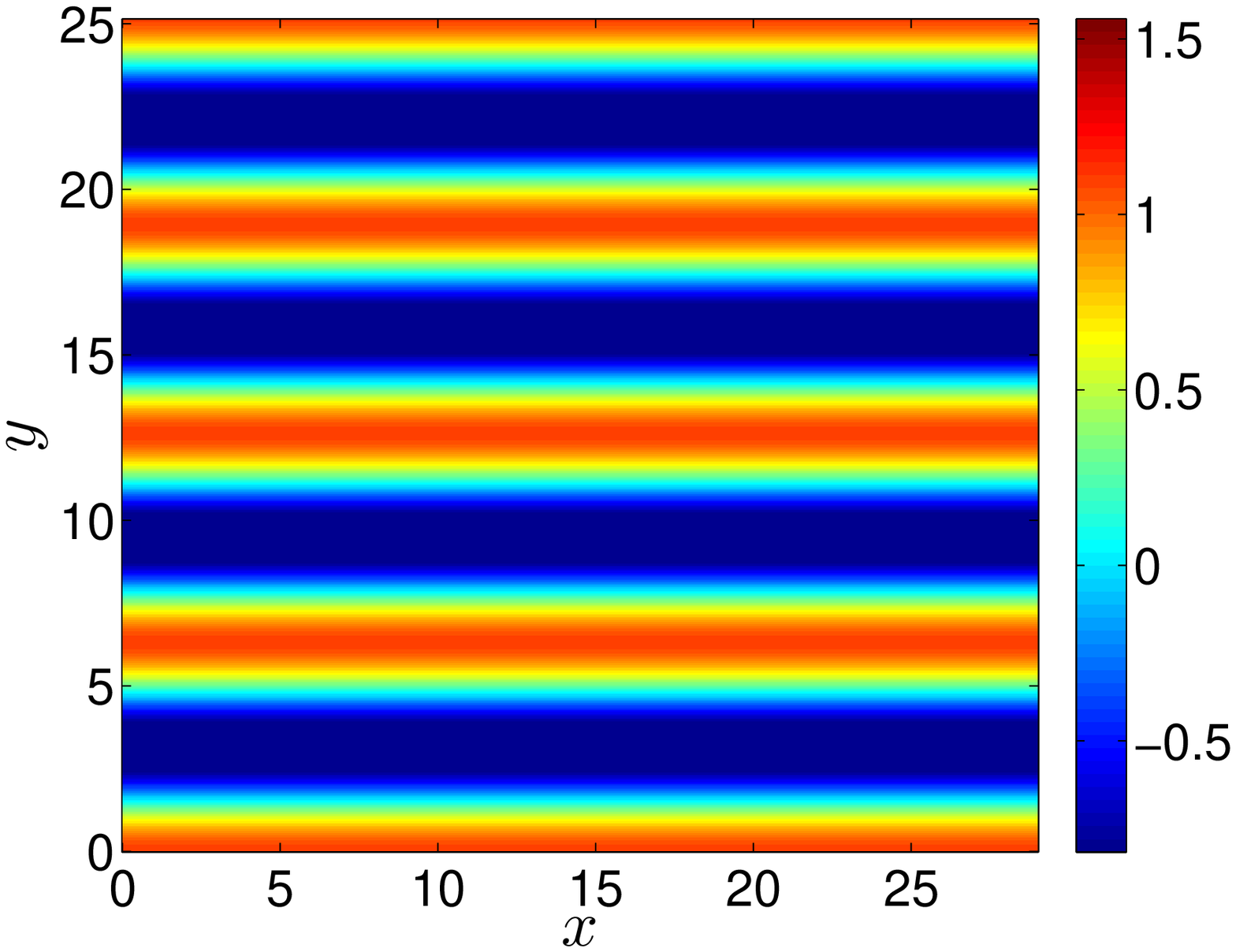}
\label{figA:2D:metastable state1}
\end{subfigure}
\hfill
\begin{subfigure}[b]{0.32\textwidth}
\includegraphics[width=\textwidth]{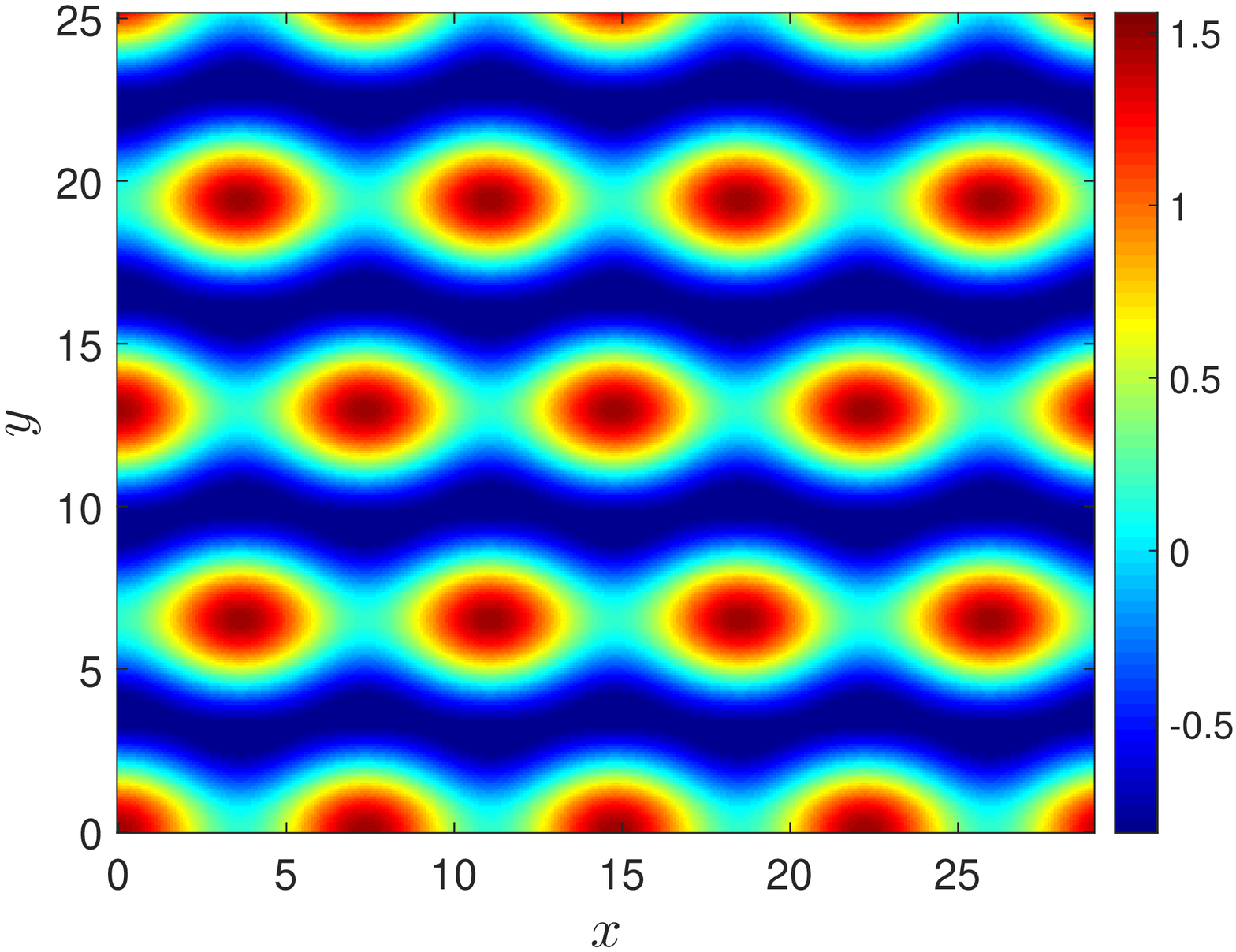}
\label{figB:2D:saddle}
\end{subfigure}
\begin{subfigure}[b]{0.32\textwidth}
\includegraphics[width=\textwidth]{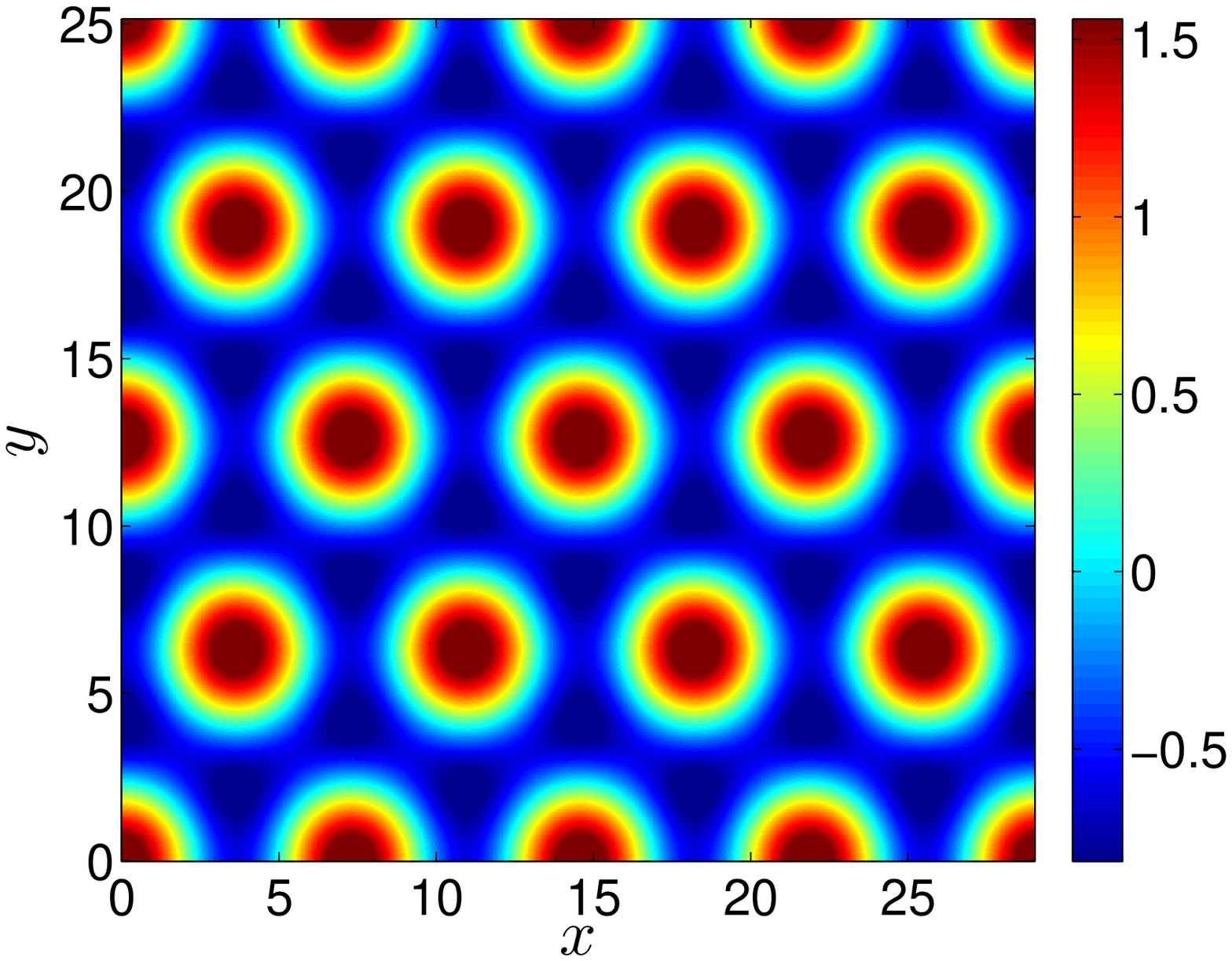}
\label{figC:2D:metastable state2}
\end{subfigure}
\hfill
\vskip -0.5cm
(a)       \hspace{3.5cm}     (b)       \hspace{3.5cm}          (c)
\caption{Two stable stationary states (a) and (c) and the transition state (b) of the 2-D  Landau-Brazovskii  energy  in $H^{-1}$ metric.
Their free energies are  $-16.486$, $-16.447$ and $-17.290$, from left to right.
Their smallest eigenvalues are $3.50\times 10^{-6}$,  $ -4.32\times 10^{-6} $, and $3.32\times 10^{-6}$, respectively.
}
\label{fig:2D_metastable_states}
\end{center}
\end{figure}

\begin{table}[htbp]
\begin{subtable}{\linewidth}\label{tab:LB_2D_inner}
\centering
{
\begin{tabular}{|*{7}{r|}}
\hline
\multicolumn{7}{|c|}{The required number of iterations} \\\cline{1-7}
{\multirow{2}*{$\Delta t$}}{}  &  \multicolumn{2}{|c|}{$err = 10^{-2}$}  & \multicolumn{2}{|c|}{$err = 10^{-3}$} & \multicolumn{2}{|c|}{$err = 2\times 10^{-4}$} \\\cline{2-7}
{} & \eqref{CSS_2D} &  \eqref{nCSS_2D} & \eqref{CSS_2D} &  \eqref{nCSS_2D} & \eqref{CSS_2D} & \eqref{nCSS_2D} \\ \hline
$0.1$ & 127 & 126 & 1930 & 1928 & 10651 & 10645 \\
\hline
$1.0$ & 14 & 14 & 192 & 190 & 1061 & 1055 \\
\hline
$10$ & 3 & 3 & 18 & 16 & 105 & 99 \\
\hline
\end{tabular}
}
\caption{ Inner iteration comparison}
\end{subtable}
\vfill
\begin{subtable}{\linewidth}\label{tab:LB_2D_outer}
\centering{
\begin{tabular}{|*{7}{r|}}
\hline
\multicolumn{7}{|c|}{The number of   cycles} \\ \cline{1-7}
{\multirow{2}*{$\Delta t$}}
  {} & \multicolumn{2}{|c|}{iter$ \#= 500 $}  & \multicolumn{2}{|c|}{iter$\#=800 $} & \multicolumn{2}{|c|}{iter$\#=1000 $} \\\cline{2-7}
{} & \eqref{CSS_2D} &  \eqref{nCSS_2D} & \eqref{CSS_2D} &  \eqref{nCSS_2D} & \eqref{CSS_2D} & \eqref{nCSS_2D} \\ \hline
 $0.1$ & 191 & 191 & 381 & 380  & 239 & 238 \\ \hline
 $1.0$ & 40 & 39 & 25 & 25 & 21 & 20 \\ \hline
 $10$ & 6 & 5 & 4 & 4 & 4 & 3 \\ \hline
\end{tabular}
}
\caption{ Outer cycle comparison }
\end{subtable}
\caption{ (a) The comparison of the CS    \eqref{CSS_2D} and  nCS scheme \eqref{nCSS_2D} for the subproblem
of the first cycle.
The integers shown are
 the  required number of iterations
 to achieve the three   prescribed    tolerances $\|\Delta\delta_\phi L(\phi^n) \|_{H^{-1}} \leq   10^{-2},
 10^{-3}$ and $  2\times 10^{-4}$; (b) The    number of  outer cycles required for \eqref{CSS_2D} and \eqref{nCSS_2D} to attain the given   tolerance $\|\Delta\delta_\phi F(\phi^{(k)})\|_{H^{-1}} \leq  10^{-6}$.  }\label{tab:LB_2D}
\end{table}

Table \ref{tab:LB_2D}   shows the number of inner iterations in the first cycle and the number of outer cycles in the whole process for the CS scheme \eqref{CSS_2D} and the nCS scheme \eqref{nCSS_2D}.
For this particular example, we find that both schemes
perform quite well for large time step sizes and their performances are almost identically.
So, in this 2D example, the both schemes are quite successful.
However, for the convex splitting scheme, we   have offered an automatic   procedure
from the decomposition of $F$ to the decomposition of $L$, so our CS  scheme developed here  inherits  the advantage of convex splitting method
to allow large time step sizes.
Note that  the linearization idea of constructing semi-implicit scheme
 has shown the  unconditionally stability   in search of local minimizers for   this Landau-Brazovskii  example.
So for other  schemes which are not based on the convex splitting idea but
have proven to  work for $F$ with unconditional stability, one may also be able to
 construct  some resemblant schemes
for $L$, such as the scheme \eqref{nCSS_2D} we derived above.
 Then what is not obvious  is the theoretic question on the stability of these new schemes for $L$.
 This could be left as a future project to generalize our idea
 of testing the convex splitting method in this article.

\subsection{Discussion} \label{ssec:dis}
The choice of the initial guess is an important practical issue
for all existing numerical methods of calculating the saddle points.
The GAD or the  IMF only has the local convergence and thus
in   extreme cases,
one can easily construct a   special initial guess which does not have convergence.
In addition, the subproblem of minimizing the auxiliary functional is well-defined
only when the minimal eigenvalue of the original Hessian is negative; otherwise,
 it should not be solved thoroughly, but limited to a fixed few number of iterations, such as
 in  the GAD. For the readers particularly interested in the practical convergence issue,
refer to the discussion in \cite{IMA2015}.
Most  numerical results we reported here are for one typical choice of initial conditions
(except for Figure \ref{CH_rate_compar}), but our unreported numerical experiments
by trying various initial guesses, still strongly support  our main conclusions.

Finally,  we emphasize  another critical implementation issue of
selecting the correct ``min-mode''.
In theory, the calculation of the min-mode is straight forward
by minimizing the corresponding Rayleigh quotient.
But if one starts from a local minimizer of $F$, say  $\bar{\phi}$,
then the eigenvalues at this locally stable state are
$\set{0=\lambda_1<\lambda_2<\ldots}$ for the periodic boundary condition and
the min-mode is then the zero eigenvector $v_1=\partial_x \bar{\phi}$.
However, taking $\partial_x \bar{\phi}$ as the min-mode
is a very bad choice since it will not push the state away
but only translate  the state    back and forth in space.
This pathological case could also appear
for some special  initial states in the convex region of $F$,
for example, when  the minimal eigenvalue $\lambda_1$ crosses over zero from positive to negative
(from the convex region to the non-convex region), i.e., near  the so-called branching point.
 The remedy to avoid this pathological situation
in practice is    simple:
at a state $\phi$,
 whenever the angle between   $v_1$ and
 $\partial_x \phi$ is  close to $0^\circ$ or $180^\circ$
(determined by a prescribed threshold),  a constraint is added to make sure that
 the min-mode $v$ in use for the auxiliary functional $L$
must be orthogonal to $\partial_x \phi$. In this way,
a strictly positive $\lambda_2$
 is selected and accordingly, $v_2$ is selected as  the ``min-mode''.
  As long as $\lambda_1$ starts to take a negative value,
  the angle defined above automatically
becomes $90^\circ$   and there is
no interference between the min-mode $v_1$ and the translation direction $\partial_x \phi$.
After taking care of this issue, we found that  for  many  initial guesses we tried,
we did observe the convergence of the algorithm.

\section{Conclusion}\label{sec5}

We have demonstrated how the convex splitting method
can improve the efficiency of the transition-state calculation
by allowing for  the preferred  large time step size.
For the 1-D Ginzburg-Landau energy,
this new method has been applied
to find index-1 saddle points
of  the Allen-Cahn and Cahn-Hilliard types, i.e., under the $L_2$ and $H^{-1}$ metrics,  respectively.  Besides, we also test this method for the 2-D Landau-Brazovskii energy functional in $H^{-1}$ metric.
The main advantage of using the convex splitting scheme
is to avoid the instability when the time step size is large.
And it is also very inspiring  that  our extensive numerical studies in this paper
have shown the significant improvement
of the computational  efficiency.
Therefore, for spatially extended systems driven by an energy functional
such as the   phase field models or  the   Kohn-Sham density functional (\cite{GAD-DFT2015}),
 we have reasons to speculate that
many matured and excellent numerical methods
for the traditional gradient dynamics
may be able to   exhibit their new vitalities
  for  saddle point calculation,  if they are correctly wrapped  by  the
  iterative minimization formulation.

\bibliography{./CVXIMF}

\begin{thebibliography}{10}

\bibitem{AllenCahn}
Samuel~M. Allen and John~W. Cahn.
\newblock A microscopic theory for antiphase boundary motion and its
  application to antiphase domain coarsening.
\newblock {\em Acta Metallurgica}, 27:1085--1095, 1979.

\bibitem{BatesFife1990}
Peter~W. Bates and Paul~C. Fife.
\newblock Spectral comparison principles for the {C}ahn-{H}illiard and
  phase-field equations and time scales for coarsening.
\newblock {\em Phys. D}, 43(2-3):335--348, 1990.

\bibitem{BatesFife1993}
Peter~W. Bates and Paul~C. Fife.
\newblock The dynamics of nucleation for the {C}ahn-{H}illiard equation.
\newblock {\em SIAM J. Appl. Math.}, 53(4):990--1008, 1993.

\bibitem{CahnHilliard}
John~W. Cahn and John~E. Hilliard.
\newblock Free energy of a nonuniform system. \textsc{I}. interfacial free
  energy.
\newblock {\em J. Chem. Phys.}, 28(2):258--267, 1958.

\bibitem{Cances2009}
E.~Canc\`{e}s, F.~Legoll, M.-C. Marinica, K.~Minoukadeh, and F.~Willaime.
\newblock Some improvements of the activation-relaxation technique method for
  finding transition pathways on potential energy surfaces.
\newblock {\em J. Chem. Phys.}, 130(11):114711, 2009.

\bibitem{cerjan1981}
C.~J. Cerjan and W.~H. Miller.
\newblock On finding transition states.
\newblock {\em J. Chem. Phys.}, 75(6):2800--2806, 1981.

\bibitem{Crippen1971}
G.~M. Crippen and H.~A. Scheraga.
\newblock Minimization of polypeptide energy : {XI}. the method of gentlest
  ascent.
\newblock {\em Arch. Biochem. Biophys.}, 144(2):462--466, 1971.

\bibitem{String2002}
W.~E, W.~Ren, and E.~Vanden-Eijnden.
\newblock String method for the study of rare events.
\newblock {\em Phys. Rev. B}, 66:052301, 2002.

\bibitem{String2007}
W.~E, W.~Ren, and E.~Vanden-Eijnden.
\newblock Simplified and improved string method for computing the minimum
  energy paths in barrier-crossing events.
\newblock {\em J. Chem. Phys.}, 126:164103, 2007.

\bibitem{GAD2011}
W.~E and X.~Zhou.
\newblock The gentlest ascent dynamics.
\newblock {\em Nonlinearity}, 24(6):1831, 2011.

\bibitem{ElseyWirth2013}
{Elsey, Matt} and {Wirth, Benedikt}.
\newblock A simple and efficient scheme for phase field crystal simulation∗.
\newblock {\em ESAIM: M2AN}, 47(5):1413--1432, 2013.

\bibitem{eyre1998unconditionally}
David~J. Eyre.
\newblock An unconditionally stable one-step scheme for gradient systems.
\newblock {\em manuscript,
  http://www.math.utah.edu/~eyre/research/methods/papers.html}, 1998.

\bibitem{eyre1998marching}
David~J. Eyre.
\newblock Unconditionally gradient stable time marching the
  \textsc{C}ahn-\textsc{H}illiard equation.
\newblock In Jeffrey~W. Bullard et~al., editors, {\em Computational and
  Mathematical Models of Microstructural Evolution}, volume 529, pages 39--46.
  MRS, Warrendale, PA,1998.

\bibitem{IMF2014}
W.~Gao, J.~Leng, and X.~Zhou.
\newblock An iterative minimization formulation for saddle point search.
\newblock {\em SIAM J. Numer. Anal.}, 53(4):1786--1805, 2015.

\bibitem{IMA2015}
W.~Gao, J.~Leng, and X.~Zhou.
\newblock Iterative minimization algorithm for efficient calculations of
  transition states.
\newblock {\em J. Comput. Phys.}, 309:69--87, 2016.

\bibitem{gu2014energy}
S.~Gu, H.~Zhang, and Z.~Zhang.
\newblock An energy-stable finite-difference scheme for the binary
  fluid-surfactant system.
\newblock {\em J. Comput. Phys.}, 270:416--431, 2014.

\bibitem{GUthesis}
Shuting Gu.
\newblock {\em On the Calculation of Transition States}.
\newblock PhD thesis, City University of Hong Kong, 2017.

\bibitem{Dimer1999}
G.~Henkelman and H.~J\'{o}nsson.
\newblock A dimer method for finding saddle points on high dimensional
  potential surfaces using only first derivatives.
\newblock {\em J. Chem. Phys.}, 111(15):7010--7022, 1999.

\bibitem{PhysRevLett.86.664}
G.~Henkelman and H.~J\'onsson.
\newblock Theoretical calculations of dissociative adsorption of
  ${\textsc{ch}}_{4}$ on an \textsc{I}r(111) surface.
\newblock {\em Phys. Rev. Lett.}, 86(4):664--667, 2001.

\bibitem{CI-NEB2000}
G.~Henkelman, B.~P. Uberuaga, and H.~J\'{o}nsson.
\newblock A climbing image nudged elastic band method for finding saddle points
  and minimum energy paths.
\newblock {\em J. Chem. Phys.}, 113(22):9901--9904, 2000.

\bibitem{HBK2005}
A.~Heyden, A.T. Bell, and F.J. Keil.
\newblock Efficient methods for finding transition states in chemical
  reactions: Comparison of improved dimer method and partitioned rational
  function optimization method.
\newblock {\em J. Chem. Phys.}, 123:224101, 2005.

\bibitem{SMWise2009}
Z.~Hu, S.M. Wise, C.~Wang, and J.S. Lowengrub.
\newblock Stable and efficient finite-difference nonlnear-multigrid schemes for
  the phase field crystal equation.
\newblock {\em J. Comp. Phys.}, 228(15):5323--5339, 2009.

\bibitem{NEB1998}
H.~J\`{o}nsson, G.~Mills, and K.~W. Jacobsen.
\newblock Nudged elasic band method for finding minimum energy paths of
  transitions.
\newblock In B.~J. Berne, G.~Ciccotti, and D.~F. Coker, editors, {\em Classical
  and Quantum Dynamics in Condensed Phase Simulations}, page 385, New Jersey,
  1998. LERICI, Villa Marigola,Proceedings of the International School of
  Physics, World Scientific.

\bibitem{KS2008}
J.~K{\"a}stner and P.~Sherwood.
\newblock Superlinearly converging dimer method for transition state search.
\newblock {\em J. Chem. Phys.}, 128:014106, 2008.

\bibitem{Kornhuber2006}
Ralf Kornhuber and Rolf Krause.
\newblock Robust multigrid methods for vector-valued {A}llen--{C}ahn equations
  with logarithmic free energy.
\newblock {\em Computing and Visualization in Science}, 9(2):103--116, Jun
  2006.

\bibitem{LOR2013}
J.~Leng, W.~Gao, C.~Shang, and Z.-P. Liu.
\newblock Efficient softest mode finding in transition states calculations.
\newblock {\em J. Chem. Phys.}, 138(9):094110, 2013.

\bibitem{GAD-DFT2015}
Chen Li, Jianfeng Lu, and Weitao Yang.
\newblock Gentlest ascent dynamics for calculating first excited state and
  exploring energy landscape of {K}ohn-{S}ham density functionals.
\newblock {\em The Journal of Chemical Physics}, 143(22):224110, 2015.

\bibitem{MMSLZZ2013}
T.~Li, P.~Zhang, and W.~Zhang.
\newblock Nucleation rate calculation for the phase transition of diblock
  copolymers under stochastic {C}ahn--{H}illiard dynamics.
\newblock {\em Multiscale Modeling \& Simulation}, 11(1):385--409, 2013.

\bibitem{Tiejun2D}
Tiejun Li, Pingwen Zhang, and Wei Zhang.
\newblock Nucleation rate calculation for the phase transition of diblock
  copolymers undr stochastic cahn-hilliard dynamics.
\newblock {\em Multiscale Model. Simul.}, 11:385--409, 2013.

\bibitem{Mousseau2000}
R.~Malek and N.~Mousseau.
\newblock Dynamics of {L}ennard-{J}ones clusters: A characterization of the
  activation-relaxation technique.
\newblock {\em Phys. Rev. E}, 22(6):7723--7728, 2000.

\bibitem{ART1998}
N.~Mousseau and G.T. Barkema.
\newblock Traveling through potential energy surfaces of disordered materials:
  the activation-relaxation technique.
\newblock {\em Phys. Rev. E}, 57:2419, 1998.

\bibitem{OKHAJ2004}
R.~A. Olsen, G.~J. Kroes, G.~Henkelman, A.~Arnaldsson, and H.~J\`{o}nsson.
\newblock Comparison of methods for finding saddle points without knowledge of
  the final states.
\newblock {\em J. Chem. Phys.}, 121(20):9776--9792, 2004.

\bibitem{Schlegel2003}
H.~Bernhard Schlegel.
\newblock Exploring potential energy surfaces for chemical reactions: An
  overview of some practical methods.
\newblock {\em J. Comput. Chem.}, 24(12):1514--1527, 2003.

\bibitem{JieShen2012}
J.~Shen, C.~Wang, X.~M. Wang, and S.M. Wise.
\newblock Second-order convex splitting schemes for gradient flows with
  {E}hrlich-{S}chwoebel type energy: application to thin film epitaxy.
\newblock {\em SIAM J. Numer. Anal.}, 50(1):105--125, 2012.

\bibitem{JieShen2014}
J.~Shen and X.F. Yang.
\newblock Decoupled energy stable schemes for phase-field model of two-phase
  complex fluids.
\newblock {\em SIAM J. SCI. Comput.}, 36(1):B122--B145, 2014.

\bibitem{CWang_mPFC2011}
C.~Wang and S.M. Wise.
\newblock An energy stable and convergent finite-difference scheme for the
  modified phase field crystal equation.
\newblock {\em SIAM J. Numer. Anal.}, 49(3):945--969, 2011.

\bibitem{RobertA}
Robert~A. Wickham, An-Chang Shi, and Zhen-Gang Wang.
\newblock Nucleation of stable cylinders from a metastable lamellar phase in a
  diblock copolymer melt.
\newblock {\em Journal of Chemical Physics}, 118:10293--10305, 2003.

\bibitem{wise2009energy}
S.~M. Wise, C.~Wang, and J.S. Lowengrub.
\newblock An energy-stable and convergent finite-difference scheme for the
  phase field crystal equation.
\newblock {\em SIAM J. Numer. Anal.}, 47(3):2269--2288, 2009.

\bibitem{SMWise2010}
S.M. Wise.
\newblock Unconditionally stable finite difference, nonlinear multigrid
  simulation of the {C}ahn-{H}illiard-{H}ele-{S}haw system of equations.
\newblock {\em J. Sci. Comput}, 44(1):38--68, 2010.

\bibitem{DuSIAM2012}
J.~Zhang and Q.~Du.
\newblock Shrinking dimer dynamics and its applications to saddle point search.
\newblock {\em SIAM J. Numer. Anal.}, 50:1899--1921, 2012.

\bibitem{ZHANG:hn}
L.~Zhang, W.~Ren, A.~Samanta, and Q.~Du.
\newblock {Recent developments in computational modelling of nucleation in
  phase transformations}.
\newblock {\em npj\ Computational Materials}, 2:16003, 2016.

\bibitem{CH1D2012}
W.~Zhang, T.~Li, and P.~Zhang.
\newblock Numerical study for the nucleation of one-dimensional stochastic
  {C}ahn-{H}illiard dynamics.
\newblock {\em Comm. Math. Sci.}, 10:1105--1132, 2012.

\end{thebibliography}

\bibliographystyle{plain} 

\end{document}